\documentclass[11pt]{article}

\title{Computational study of irrational rotations via exact discontinuity tracking: A new experimental framework%
\author{Hannah Kravitz \footnote{hkravitz@pdx.edu} \\
{Published in \textit{Experimental Mathematics},
Oct 2026}. \\
{DOI: \href{https://doi.org/10.1080/10586458.2026.2742231}
{10.1080/10586458.2026.2742231}.}}}
\date{}

\usepackage{geometry}
\usepackage{amsmath, amssymb} 
\usepackage{graphicx}         
\usepackage{float}            
\usepackage{algorithm}        
\usepackage{algpseudocode}    
\usepackage{subfigure}
\usepackage[noadjust]{cite}
\usepackage{url}
\usepackage{xcolor}
\usepackage[colorlinks=true, allcolors=blue]{hyperref}
\usepackage{epstopdf}

\begin{document}

\maketitle

\begin{abstract}
\textcolor{black}{We present an algorithm for studying the sum $$ G_N(x,\rho) = \sum_{k=1}^N \left[ \{x+k\rho\} - \frac{1}{2}\right],$$ which is closely related to the discrepancy of the irrational rotation. Computation of $ G_N $ and associated statistics is notoriously difficult due to numerical instability and the failure of sampling to capture the discontinuities. Our algorithm instead fully defines $G_N$ through its discontinuities, which allows computation of $ G_N (x,\rho)$ with errors on the order of $N$ times machine epsilon in $O(N)$ time for practical $N$. The associated probability density function can then be computed in $O(N \log N)$ time.}

\textcolor{black}{
Presented as a tool for experimental mathematics, the algorithm enables reliable plots of the descriptive statistics of $G_N$ over $N$, in particular the length of the support of the $\text{pdf}$ (the discrepancy), revealing extremal behavior related to the arithmetic properties of $\rho$. We also observe a new pattern: when $\rho$ is well approximated by a continued-fraction convergent $\frac{p_n}{q_n}$, the $\text{pdf}$ exhibits a predictable spiked-trapezoidal pattern when $N=kq_n$, eventually degrading at a rate that appears to depend on approximation quality.}
\end{abstract}


\section{Introduction}
\label{sec:intro}

\textcolor{black}{The key object of interest in this paper is the centered sum of fractional parts}
\begin{equation}
\label{eq:discrepancy_sum}
\color{black}G_N\color{black}(x,\rho) = \sum_{k=1}^N \left[ \color{black}\{x+k\rho\}\color{black} - \frac{1}{2}\right].
\end{equation}
with $x \in [0,1)$, $\rho$ (typically irrational) $\in [0,1)$, and $N \in \mathbb{N}$, \textcolor{black}{and where $\{y\}= y - \lfloor y \rfloor$. $G_N$ appears across many different mathematical settings. It} is a particular example of a Birkhoff sum\textcolor{black}{; when $\rho$ is irrational, $\frac{G_N}{N}$} quantifies the deviations in the average of the sequence $x_k=\color{black}\{ x+k\rho \}$ from its expected mean of $\frac{1}{2}$~\cite{birkhoff1931proof}. 
\textcolor{black}{Another interpretation of $G_N$ comes from quasi-Monte Carlo integration~\cite{niederreiter1978quasi}; $\frac{G_N}{N}$ quantifies the difference between the equal-weight quadrature 
$\frac{1}{N} \sum_{k=1}^N x_k$
with nodes $x_k=\{x+k \rho \}$ and the integral 
$\int_0^1 t \,dt = \frac{1}{2}$. When $\rho$ is irrational, $\{x+k \rho \}$ is a classic example of a uniformly distributed modulo one sequence~\cite{kuipers1974uniform}, which implies
\begin{align*}
       \lim_{N \to \infty} \frac{1}{N} \sum_{k=1}^N \{x+k \rho \} = \int_0^1 t \,dt = \frac{1}{2}
\end{align*}
i.e. 
\begin{align*}
\lim_{N \rightarrow \infty} \frac{G_N(x,\rho)}{N} =0 \text{ for } \rho \notin \mathbb{Q}.
\end{align*}
$G_N$ is also closely related to the discrepancy of the irrational rotation by the identity~\cite{ramshaw1981discrepancy, baxa1994minimum, pinner1997sums}
\begin{equation} \label{eq:disc_id}
2 \|G_N (\cdot, \rho)\|_{\infty}  =D_N \left( \rho  \right)
\end{equation}
where
$$D_N \left( \rho \right)= \sup_{[a,b) \subseteq [0,1)} \left| \sum_{k=1}^N \mathbf{1}_{[a,b)} \left( \{k\rho \}\right) -N (b-a) \right|$$
where $\mathbf{1}_{[a,b)}$ is the indicator function of the half-open interval $[a,b)$.}

\textcolor{black}{These connections place $G_N$ within a substantial theoretical literature, with} work spanning nearly a century \textcolor{black}{analyzing} growth rates, extrema, asymptotic behavior, \textcolor{black}{bounds}, and other properties of $\color{black} G_N \color{black}(x,\rho)$ \textcolor{black}{and closely related objects} (see, for example,~\cite{antonevich2022behaviour, bountis2020cauchy, doi2017upper, knill2011self, kochergin2023growth, kuipers1974uniform, mori2019distribution, ramshaw1981discrepancy, roccadas2011local, schoissengeier1986discrepancy, setokuchi2014discrepancies, setokuchi2015discrepancy, shimaru2017discrepancies, shimaru2018discrepancies, das2016quasiperiodicity, tong2026weighted, pinner1997sums, krasil2011description, baxa1994minimum, borda2023limit}) yet many questions still remain unanswered. 
A primary reason for this enduring interest is that this simple function exhibits highly complex self-similar patterns that evolve and change, determined by the \textcolor{black}{arithmetic} properties of $\rho$ and the choice of $N$. However, these features grow finer as $N$ increases, making them difficult to resolve computationally using standard methods. This creates a significant barrier to using numerical exploration for testing conjectures. Furthermore, the theoretical literature on this topic is often presented without figures, obscuring the fascinating geometric patterns the mathematics reveal.

The computational challenges are threefold:
\begin{itemize}
    \item First, the naive summation in Equation~\ref{eq:discrepancy_sum} is notoriously prone to numerical instability; the repeated use of the $\text{mod }{1}$ operation on large terms runs a high risk of catastrophic cancellation (loss of nearly all significant digits from subtracting two nearby numbers)~\cite{cuyt2001remarkable, heath2018scientific}.
    \item Second, the function's structure consists of up to $N\color{black}+1$ branches with jump discontinuities. For irrational $\rho$, a uniform sampling of $x \in[0,1)$ is guaranteed to miss the branch endpoints, which are fundamental to the structure of the probability density function ($\text{pdf}$). As $N$ varies, $\color{black} G_N \color{black}$ presents intricate fine-scale features for irrational $\rho$ that sampling is also nearly guaranteed to miss, like very short branches and branches that are nearly (but never exactly) aligned. 
    \item Third, resolving this complex structure to even roughly discern some features is computationally prohibitive itself, requiring \textcolor{black}{at least quadratic time and dense sampling}. This makes even familiar rotation numbers like $\rho =\pi-3$ difficult to study, as the most interesting behavior typically occurs when $N$ is a multiple of \textcolor{black}{the denominator of a continued fraction convergent; the fourth such convergent denominator} for $\pi-3$ is $N=33,102$.
\end{itemize}

In this paper, we overcome all three of these barriers by introducing an efficient and stable algorithm that computes $\color{black} G_N \color{black}$ to machine precision. By developing a novel method to track the $N$ discontinuities directly, our algorithm achieves what was unattainable using naive sampling: it runs in \textcolor{black}{linear time for practical $N$}, operates with machine precision (bypassing the numerical instabilities inherent in the summation and avoiding large sampling errors) \textcolor{black}{with errors on the order of $N$ times machine epsilon ($\varepsilon_{\text{mach}}$)}, and requires only $O(N)$ storage.

A direct consequence of this algorithm is the ability to exactly compute the probability density function ($\text{pdf}$) associated with $\color{black} G_N \color{black}(x,\rho)$ in $O(N \log N)$ time, with $O(N)$ storage, at machine precision. This new computation enables the exact calculation of the fundamental statistical properties of the $\text{pdf}$: $||\color{black} G_N \color{black}||_{\infty}$ (half of the \textcolor{black}{length of the} support of the $\text{pdf}$\textcolor{black}{, or equivalently half of $D_N$}), $||\color{black} G_N \color{black}||_{2}^2$ (the variance of $\text{pdf}$), and the kurtosis. \textcolor{black}{To our knowledge, this is the first algorithm to compute the discrepancy $D_N(\rho)$ directly from the fully resolved spatial structure of $G_N(x,\rho)$. This provides some numerical advantages, namely controlled, predictable error. Resolving the spatial structure using sampling would require extremely dense sampling, so this also ensures a quick runtime and minimal storage.}

In the context of experimental mathematics, the ability to visualize a mathematical object is often a precursor to understanding it. This algorithm provides the necessary stability and resolution to reveal subtle \textcolor{black}{geometric} structures of the \textcolor{black}{sum} and its \textcolor{black}{statistical moments} that would remain obscured by standard numerical measures.

\subsection{\textcolor{black}{Structure of $G_N$}}

Equation~\ref{eq:discrepancy_sum} may be rewritten in a form more stable to double-precision arithmetic errors:
\begin{equation}  \label{eq:disc_sum_rewrite}
\color{black} G_N \color{black}(x,\rho) = Nx + \frac{N(N+1)}{2} \rho  - \frac{N}{2} - \sum_{k=1}^N \lfloor x + k \rho \rfloor 
\end{equation}
where $\lfloor \cdot \rfloor$ is the floor function.  This formulation also reveals the structure of $\color{black} G_N \color{black}$ as piecewise-linear with branches of slope $N$. For irrational $\rho$, the jumps between segments caused by the floor function will have height exactly 1, while the lengths of the segments \textcolor{black}{vary according to the Steinhaus conjecture (Three-gap theorem)~\cite{van1988three, polanco2016continuous}, i.e. they can take on at most three values -- when three distinct gap lengths occur, the largest is the sum of the smaller two}. As $N$ increases, the plots of $\color{black} G_N \color{black}$ reveal a self-similar structure\footnote{\textcolor{black}{Here, the phrase ``self-similar'' is used heuristically to describe patterns that recur across arithmetic scales associated with convergent denominators and gap lengths of $\rho$. It is not used as a rigorous mathematical definition, and no particular scaling law is asserted.}} with intricacies determined by the \textcolor{black}{arithmetic} properties of $\rho$.

\textcolor{black}{To describe these arithmetic properties, we use standard continued fraction notation. For irrational $\rho \in (0,1)$, we write
$$\rho = [0; a_1, a_2, a_3, \dots ]$$
with $a_n \in \mathbb{N}=\{1, 2, 3, \dots \}$. The $n$th convergent is obtained by truncating this expansion:
$$\frac{p_n}{q_n} = [0; a_1, a_2, a_3, \dots a_n], \quad \gcd(p_n, q_n)=1.$$
Here, $a_n$ is the $n$th partial quotient and $q_n$ is the $n$th convergent denominator. Each convergent $\frac{p_n}{q_n}$ approximates $\rho$ more closely than any rational number with a smaller positive denominator.}

\begin{figure}[H]
    \centering
    \includegraphics[width=1\linewidth]{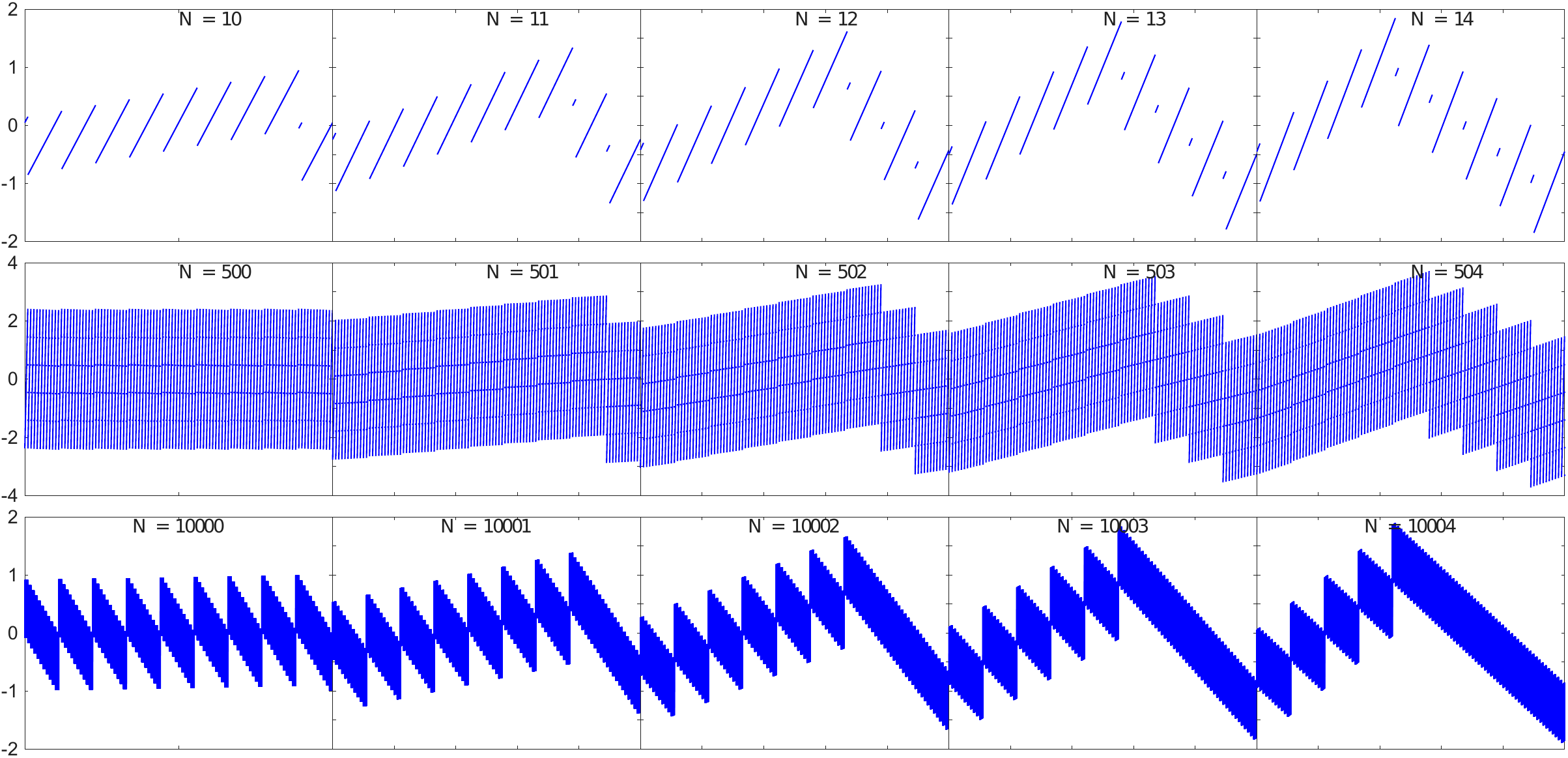}
    \caption{This figure shows the plot of $\color{black} G_N \color{black}(x,L)$ for different $N$. The three rows show different levels of a self-similar structure changing with $N$.} \label{fig:branches}
\end{figure}

Figure~\ref{fig:branches} shows $\color{black} G_N \color{black}$ for $\displaystyle \rho=L = \sum_{k=1}^\infty \frac{1}{10^{k!}} \approx 0.1100010000000000$, Liouville's constant for several $N$ to show both the capabilities of the algorithm and the behavior of the function. \textcolor{black}{$\rho=L$ provides a particularly good test case for the algorithm, as its non-Diophantine character~\cite{Khinchin_CF} can produce extremely small gap lengths~\cite{van1988three, polanco2016continuous} and therefore short branches of $G_N$.} As $N$ increases in the second and third rows, a self-similar structure emerges, with branches grouped in blocks with \textcolor{black}{repeating structure}. 

\subsection{Probability density function and its statistical properties}
\label{sec:statistics}

We find the probability density function $\text{pdf}(y)$ using the standard formula for a non-monotonic transformation on a uniform distribution~\cite{pitman2012probability}, in this case $\color{black} G_N \color{black}(x,\rho):[0,1) \rightarrow \mathbb{R}$, with $\rho$ fixed:
\begin{equation}
\text{pdf}(y)= \sum_{x \in \color{black} G_N \color{black}^{-1}(y)} \frac{1}{|\color{black} G_N \color{black}'(x, \rho)|} = \sum_{x \in \color{black} G_N \color{black}^{-1}(y)} \frac{1}{N}  = \frac{\# \color{black} G_N \color{black}^{-1}(y) }{N}
\end{equation}
where $\color{black} G_N \color{black}^{-1}(y)$ is the preimage of $y$ and $\#$ is the cardinality of a set.
The structure of $\text{pdf}(y)$ is piecewise-constant, with segments becoming extremely short as $N$ grows large.

\begin{figure}[H]
\centering
\subfigure{

 \includegraphics[angle=-90, width=0.44\linewidth]{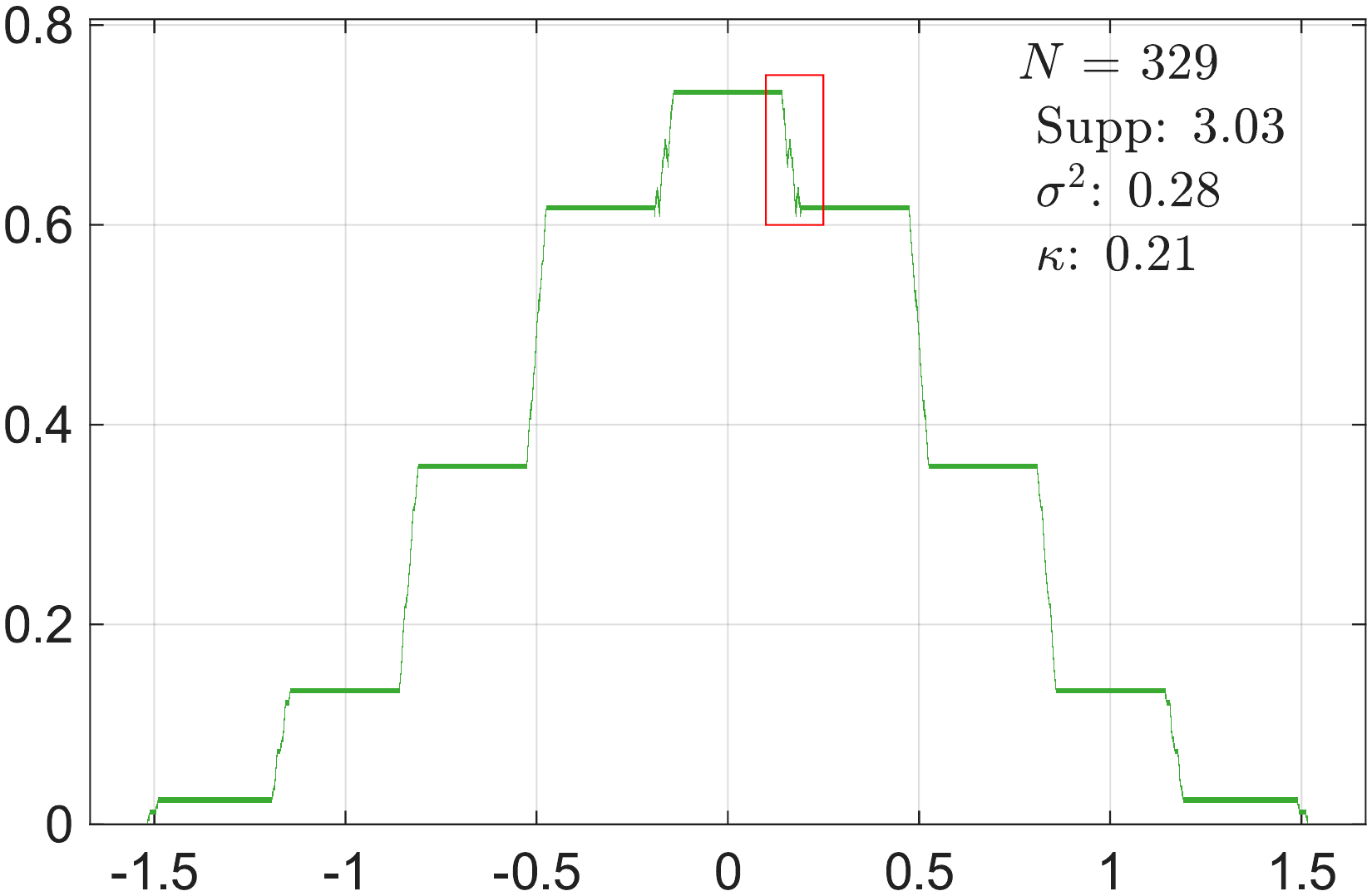}}%
\hfill
\subfigure{
    \includegraphics[angle=-90, width=0.44\linewidth]{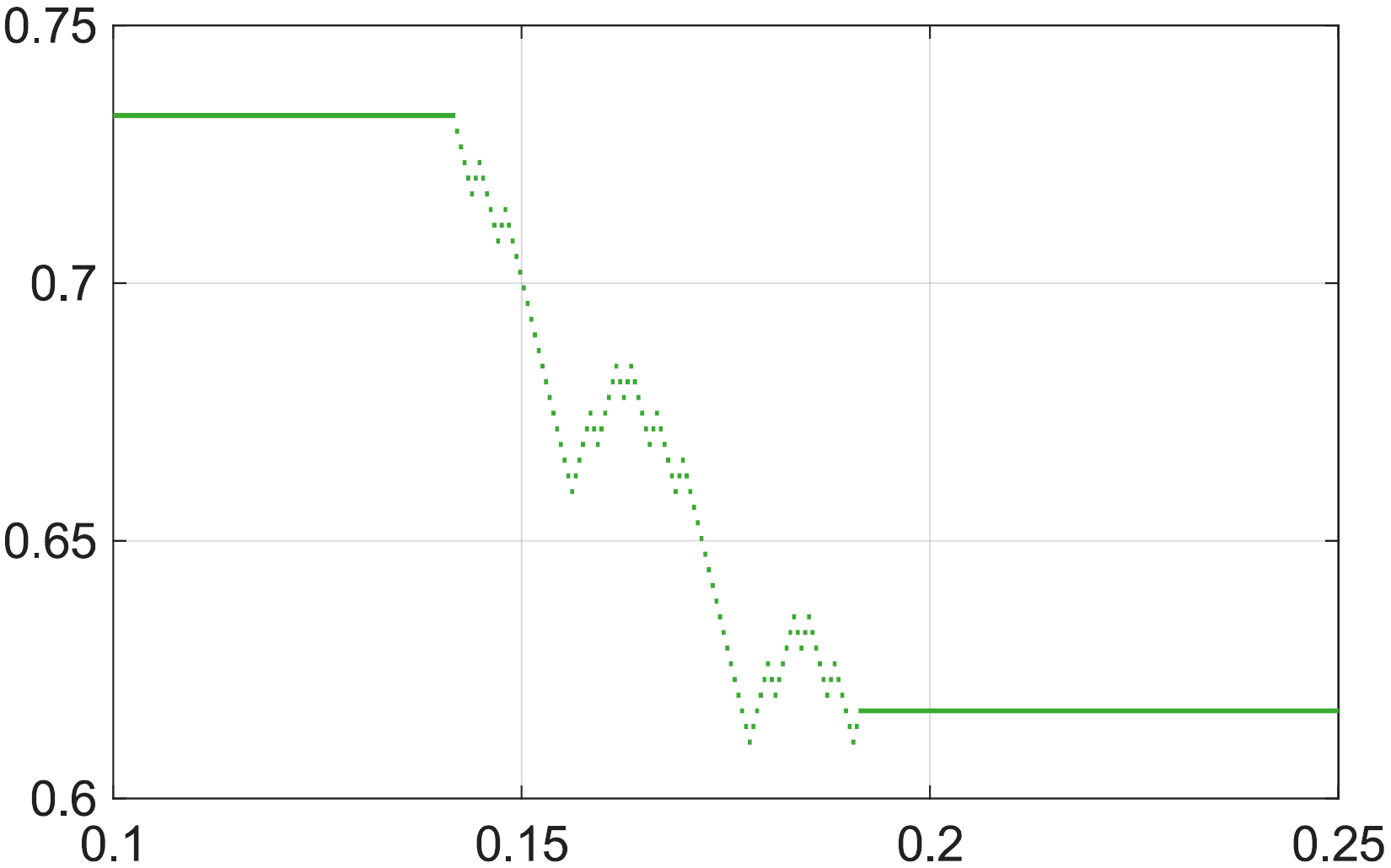} } \\
    \subfigure{
    \includegraphics[angle=-90, width=0.44\linewidth]{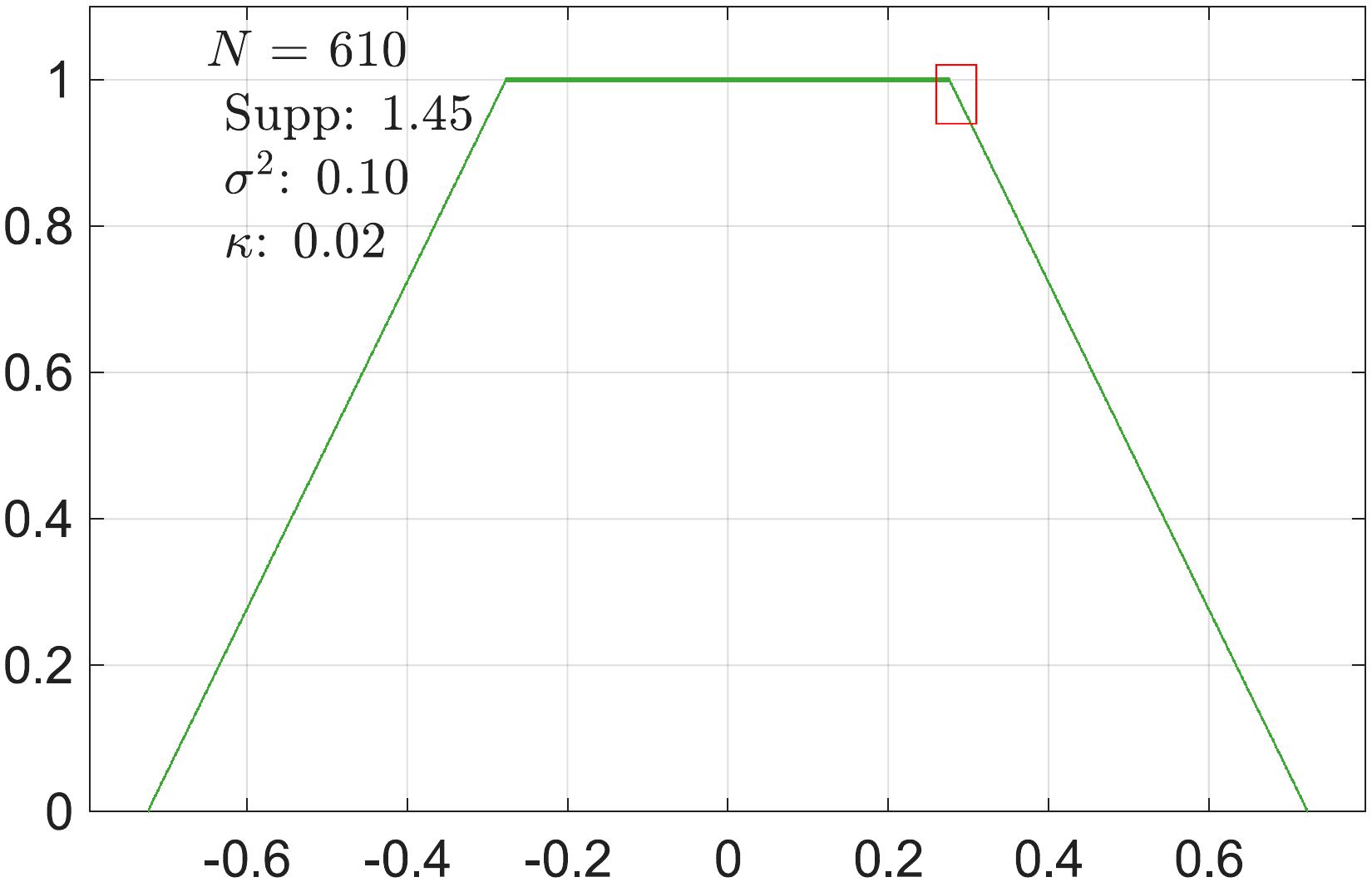}}%
\hfill
\subfigure{
    \includegraphics[angle=-90, width=0.44\linewidth]{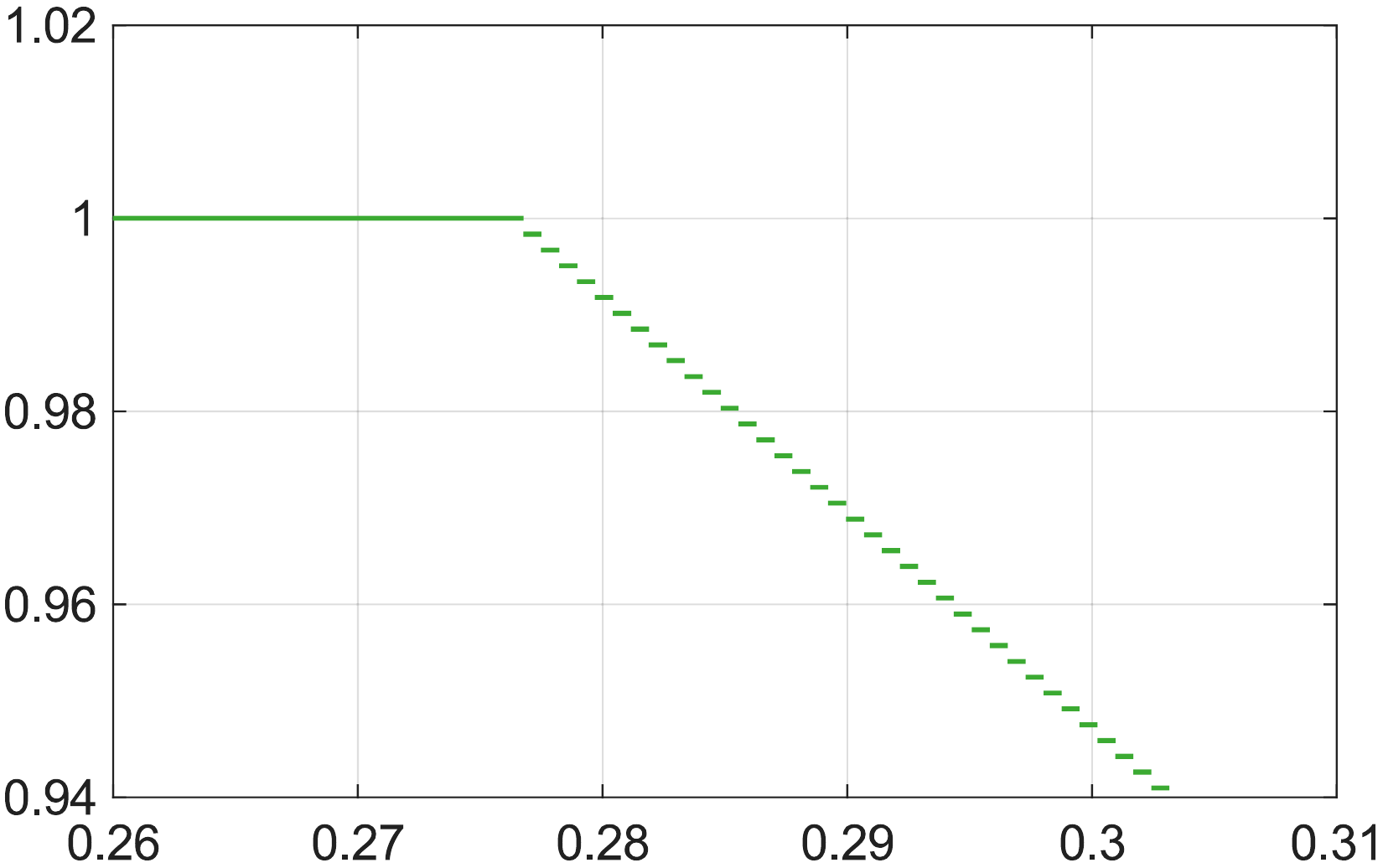} }
\caption{(Left) $\text{pdf}$ for $\rho= \frac{\sqrt{5}-1}{2}$ for two choices of $N$. The top, $N=329$, shows a few straight platforms connected by very short segments of varying heights, while $N=610$ (a \textcolor{black}{convergent denominator} for $\rho$) has a long straight horizontal line at $\color{black}\text{pdf}(y)\color{black}=1$, indicating that $\color{black} N$ branches pass through a particular range of values containing zero, with piecewise-constant line segments making up the sides of this trapezoid-like distribution. (Right) The $\text{pdf}$ functions are zoomed in to show the detail the algorithm is able to capture.}
\label{fig:pdf_rho}
\end{figure}

\begin{figure}[H]
\centering
\subfigure{
    \includegraphics[angle=-90,width=0.44\linewidth]{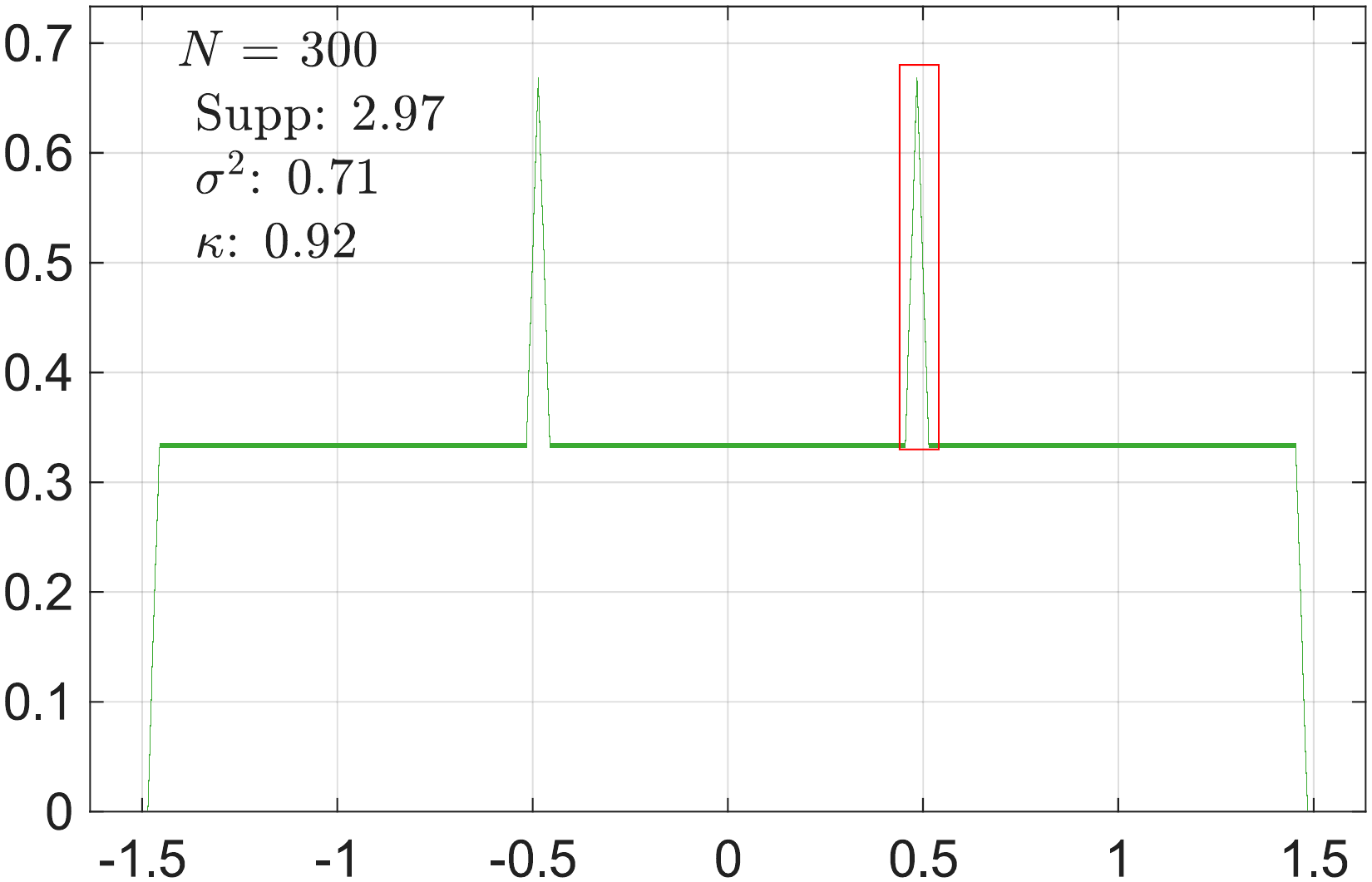}}%
\hfill
\subfigure{
    \includegraphics[angle=-90,width=0.44\linewidth]{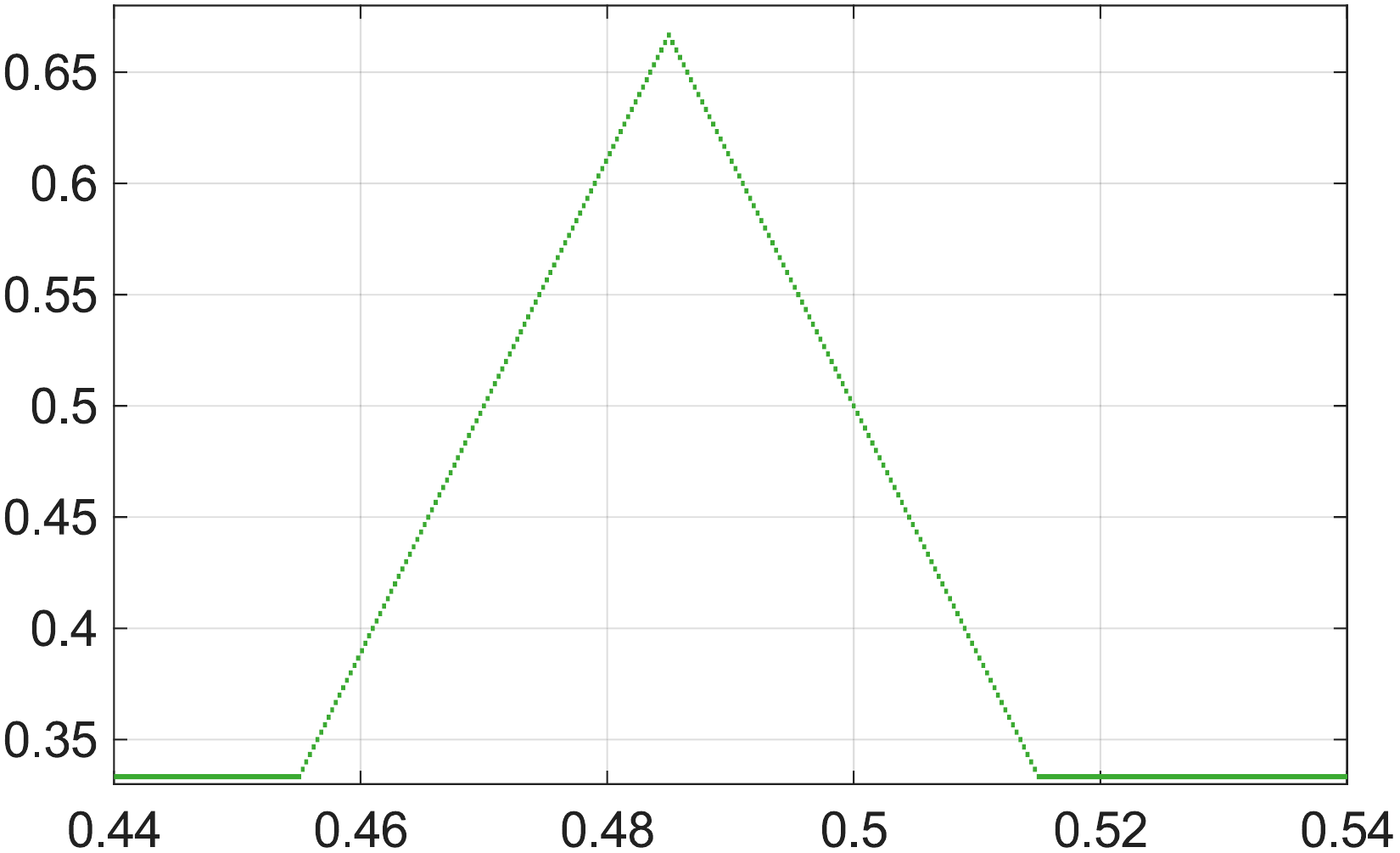} } \\
    \subfigure{
    \includegraphics[angle=-90,width=0.44\linewidth]{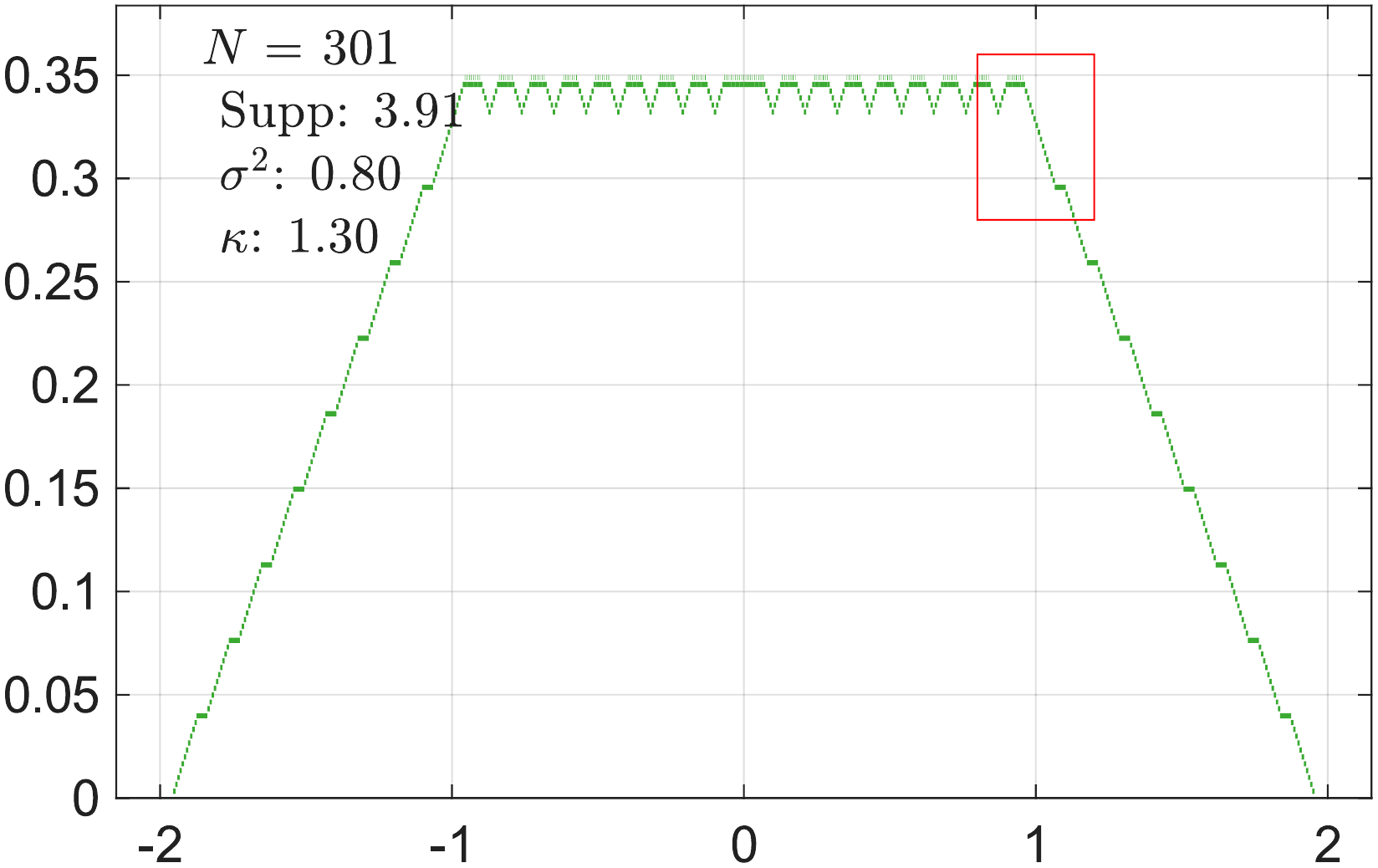}}%
\hfill
\subfigure{
    \includegraphics[angle=-90,width=0.44\linewidth]{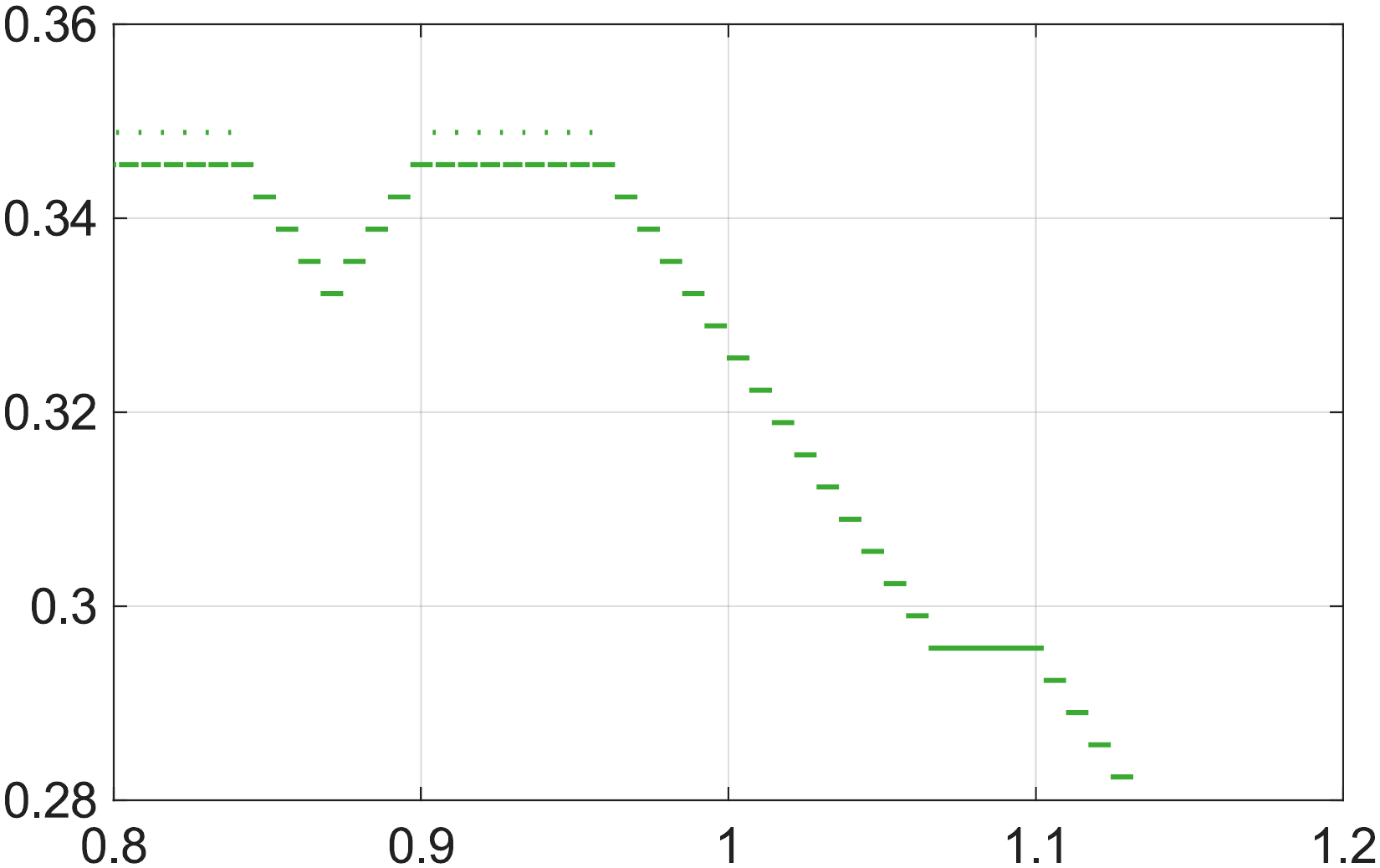} }
\caption{(Left) $\text{pdf}$ for $\rho =L$ (Liouville's constant). Two consecutive values of $N$ are shown, $N=300$, \textcolor{black}{three times} the second convergent denominator \textcolor{black}{$q_2=100$}, and $N=301 \color{black} = 3 q_2+1$, where the support of the distribution immediately spreads out. (Right) The $\text{pdf}$ functions are zoomed in to show the detail the algorithm is able to capture.}
\label{fig:pdf_Liou}
\end{figure}

While this distribution is often studied in the context of ergodic theory, \textcolor{black}{e.g.} ``almost everywhere'' and asymptotic results~\cite{dolgopyat2018no, dolgopyat2017temporal, dolgopyat2020quenched, bromberg2018temporal, fruhwirth2024birkhoff, bountis2020cauchy}, highly structured and nontrivial behaviors occur for particular $\rho$ and $N$ combinations. These patterns, which our algorithm computes with machine precision, are \textcolor{black}{related to} the \textcolor{black}{arithmetic} properties of $\rho$. A few different pattern types are shown in Figures~\ref{fig:pdf_rho} and~\ref{fig:pdf_Liou}. The left panels of the two figures show the distributions, while the right panels show 
zoomed-in versions, showing the detail that \textcolor{black}{uniform sampling generically misses}.

\subsubsection{Descriptive statistics} \label{sec:desc_stats}
We now introduce the descriptive statistics of $\text{pdf}(y)$ that can now be computed directly, while also noting their interpretation in terms of $\color{black} G_N \color{black}$ \textcolor{black}{and the discrepancy}. The descriptive statistics for particular $\text{pdf}$s are indicated in Figures~\ref{fig:pdf_rho} and~\ref{fig:pdf_Liou}.

\begin{itemize}
    \item \textbf{Support}. The support of $\text{pdf}(y)$ is the interval $\left[ \min(\color{black} G_N \color{black}), \sup(\color{black} G_N \color{black}) \right]$. \textcolor{black}{We label the length of the support $\operatorname{Supp}(\text{pdf})$.}    
    Due to the symmetry of the distribution,    
    $$\operatorname{Supp}(\text{pdf})=2\| \color{black} G_N \color{black} \|_{\infty} = 2\sup_x \color{black} G_N \color{black}(x,\rho) = - 2\min_x \color{black} G_N \color{black}(x,\rho) \color{black}.$$    
    Our algorithm also directly computes the previously-inaccessible 
    $x$-locations of the jump discontinuities \textcolor{black}{associated with these extrema}. \textcolor{black}{$\operatorname{Supp}(\text{pdf})$ also directly gives us the discrepancy (see Equation~\ref{eq:disc_id}).}
    \item \textbf{Variance}, $\sigma^2$, given by $E[Y^2]-E[Y]^2$ where $E$ is the expected value~\cite{casella2024statistical}. Since the mean of our distribution is zero, we apply the definition of the variance of a random variable $Y$ and then the Law of the Unconscious Statistician~\cite{casella2024statistical}, noting that $[0,1)$ has a uniform probability density function:
    \begin{equation} \label{eq:variance}
    \sigma^2 = E[Y^2] = \int_{-\infty}^{\infty} y^2 \text{pdf}(y) \,dy = \int_0^1 \left( \color{black} G_N \color{black}(x,\rho) \right)^2 \cdot 1 \,dx = ||\color{black} G_N \color{black}||_2^2.
    \end{equation}
    \item \textbf{Kurtosis}, $\kappa$, is a measure of ``tailedness'' of a \textcolor{black}{distribution}, with a Gaussian having kurtosis of 3~\cite{casella2024statistical}:
    \begin{equation} \label{eq:kurtosis}
        \kappa = \frac{E[Y^4]}{\sigma^4} = \frac{1}{\sigma^4}  \int_{-\infty}^{\infty}y^4 \text{pdf}(y)\,dy.
    \end{equation}
\end{itemize}

These statistics as a function of $N$ are shown in Appendix~\ref{apx:stats_plots} for several $\rho$.

\section{Numerical implementation} \label{sec:numerical_implementaion}
Several authors have studied $\color{black} G_N \color{black}$, or related objects, numerically, typically focusing on particular quantities or specific subsequences of $N$. For example, Setokuchi \textcolor{black}{computed the star discrepancy --- closely related to $G_N(0,\rho)$~\cite{baxa1994minimum} --- for $N$ of order $10^7$ by sampling every $10^4$ points~\cite{setokuchi2015discrepancy}}, while Mori, Shimaru, \& Takashima sampled multiples of \textcolor{black}{convergent denominators}~\cite{mori2019distribution} \textcolor{black}{in their study of $G_N(0,\rho)$}. An investigation of \textcolor{black}{the distributional properties of} $\color{black} G_N \color{black}(x,\rho)$ for particular $N$ and $\rho$ was conducted by Bountis, Veerman, \& Vivaldi by sampling $x$~\cite{bountis2020cauchy}, and a similar sampling-based strategy was employed by Vinson in the study of a related function~\cite{vinson2001partial}. Methods based on sampling the domain tend to be very slow and present difficulties in distinguishing mathematical behavior from numerical artifacts (see Section \ref{sec:naive}). Other numerical \textcolor{black}{and theoretical} work considers weighted\textcolor{black}{, often normalized,} variations of Equation~\ref{eq:discrepancy_sum} \textcolor{black}{and related objects}, often in the Fourier domain~\cite{blessing2024weighted, das2016quasiperiodicity, fukuyama2016metric, fukuyama2023metric, gonzalez2022efficient, tong2024exponential, das2018super, tong2026weighted}, \textcolor{black}{sometimes achieving polynomial or even exponential convergence}. In contrast, our focus is resolving the fine spatial structure of $\color{black} G_N \color{black}$ itself, enabling $\color{black} G_N \color{black}$, its $\text{pdf}$, and the associated statistics to be studied with exact, rather than sampled, information.

\subsection{Naive algorithm} \label{sec:naive}
To establish a baseline for comparison, we first consider the naive algorithm (Algorithm~\ref{alg:naive}), which directly implements Equation~\ref{eq:discrepancy_sum} by sampling $N_{\text{samples}}$ points across an interval $x \in [0,1)$. 

\begin{algorithm}[H]
\caption{Naive algorithm for $\color{black} G_N \color{black}$}\label{alg:naive}
\begin{algorithmic}
\Require $N \in \mathbb{N}$, $\rho \in [0,1)$, \texttt{Nsamples} $\in \mathbb{N}$ \Comment{\texttt{Nsamples} $\gg N$ to resolve some structure}
\State xlist $\gets$ linspace$(0,1,$Nsamples$)$
\State Slist $\gets$ zeros$(1,$Nsamples$)$

\For{$k \gets 1:$len$($xlist$)$}
    \State $S \gets0$ 
	\For{$i \gets 1:N$}
		\State $S \gets S +$ mod$($xlist$[k]+i \rho, 1) - \frac{1}{2}$
	
	\EndFor
\State Slist$[k] \gets S$

\EndFor
\end{algorithmic}
\end{algorithm}

As discussed in Section~\ref{sec:intro}, the repeated modular arithmetic produces imprecise results. \textcolor{black}{The number of samples should be much greater than the number of branches to have any chance of resolving the structure, so the storage requirement is $O(N_{\text{samples}})$ where $N_{\text{samples}} \gg N$, with time complexity $O(N N_{\text{samples}})$. In addition, generic uniform sampling is based on rational numbers while the discontinuities occur at irrational points for irrational $\rho$. Random sampling would also generically miss the key features}. 


\subsection{Discontinuity-tracking algorithm} \label{sec:DTA}
To overcome the fundamental limitations of sampling, we present a discontinuity-tracking algorithm (we use DTA for short) to exactly compute $\color{black} G_N \color{black}(x,\rho)$ with errors on the order of $\color{black} N \varepsilon_{\text{mach}}$. We first consider the structure of the equivalent expression for $\color{black} G_N \color{black}$, reproduced here for convenience:
\begin{equation*} 
\color{black} G_N \color{black}(x,\rho) = Nx + \frac{N(N+1)}{2} \rho  - \frac{N}{2} - \sum_{k=1}^N \lfloor x + k \rho \rfloor.
\end{equation*}

Where the naive implementation of modular arithmetic suffers from catastrophic cancellation, the floor function $\lfloor y \rfloor$ results in an exact integer. The expression $\frac{N(N+1)}{2}$, which multiplies $\rho$ is also an exact integer for practical $N$ (that is, $N<2^{\color{black}27}$); multiplying by $\rho$ last keeps the arithmetic in the integer regime for as long as possible. 

The main innovation in the DTA is in noticing that $\color{black} G_N \color{black}$ is piecewise-linear with branches of slope $N$ and therefore fully defined by its branch endpoints; thus we significantly reduce both the computational and storage load by computing only those points. Jump discontinuities occur at $x=\lceil k \rho \rceil - k \rho$, where $\lceil \cdot \rceil$ represents the ceiling function. This computation for $x$ may introduce some catastrophic cancellation, but it is important to note that these errors do not accumulate, as each $x$ is independent of the next. 

Another important insight is noticing that the jumps must always be of integer height, and of height 1 for irrational $\rho$. The algorithm identifies and stores discontinuities within $[0,1)$, preserving duplicates in the rational case to track the height of the jumps. First, the $y$-intercept, $\color{black} G_N \color{black}(0,\rho)$, is computed. Then, the differences between successive discontinuities ($\Delta x$) determine branch widths along the $x$-axis. $r$ duplicate discontinuities indicate a jump of height $r$. $\color{black} G_N \color{black}$ can then be fully characterized by leveraging the slopes $N$ and branch widths $\Delta x$. Minimal storage is required, as the only values that need to be stored are the branch endpoints. Algorithm~\ref{alg:main} provides a detailed implementation of this approach.

\begin{algorithm}[H]
    \caption{Discontinuity-tracking algorithm} \label{alg:main}
    \begin{algorithmic}[1]
    \Require $N \in \mathbb{N}, \rho \in [0,1)$
    \color{black} \State $k \gets 1:N$ 
    \State $\texttt{disc\_list} \gets \lceil k\rho\rceil-k\rho$ \Comment{Identify locations of discontinuities, preserving duplicates}  \color{black} 
     \color{black} \State $\texttt{disc\_list} \gets \texttt{disc\_list}[\texttt{disc\_list}>0]$
    \Comment{Remove zeros}  \color{black} 
    \State $\texttt{disc\_list} \gets $ sort(\texttt{disc\_list}, ascending)
    \State  $\color{black} \texttt{disc\_list} \gets\left[ 0, \texttt{disc\_list}\right] $ \Comment{\color{black} Prepend to sorted list \color{black}}
    \State Create $\Delta x_{\texttt{list}} \gets$ diff $([$\texttt{disc\_list},\color{black} 1$])$ \Comment{Temporarily appending the endpoint of 1 allows us to count the final width} \color{black} 
    \State $ \color{black} G_N\color{black}(0,\rho) \gets \frac{N(N+1)}{2} \rho - \frac{N}{2} - \sum_{k=1}^{N} \lfloor k \rho \rfloor$ \Comment{Compute $y$-intercept}
    \State $a_1 \gets \color{black} G_N\color{black}(0,\rho)$ and $b_1 \gets a_1 + N \Delta x_{\text{list}}[1]$. \Comment{1st branch spans $y \in [a_1, b_1)$}
    \For{$j = 2 \text{ to }$ len(\texttt{disc\_list})}
        \State Set \(a_j \gets b_{j-1} - 1\) and \(b_j \gets a_j + N \Delta x_{\text{list}}[j]\). \Comment{$\color{black}j$th branch spans $y \in [a_{\color{black}j}, b_{\color{black}j})$}
\EndFor
    \end{algorithmic}
\end{algorithm}

\textcolor{black}{By the Steinhaus conjecture (Three-gap theorem)~\cite{van1988three, polanco2016continuous, ravenstein1990characteristics}, the values in $\Delta x_{\texttt{list}}$ can take on at most three distinct values. Further, if the continued fraction expansion of $\rho$ is known, the gap widths~\cite{polanco2016continuous} and their order ~\cite{ravenstein1990characteristics} may be computed directly rather than subtracting neighboring discontinuity locations. Here we leave the algorithm generic (not requiring computation of the continued fraction structure of $\rho$) and save these improvements for future work.}

To quantify the efficiency of the improved algorithm, we first determine the number of floating point operations (FLOPs) required. \textcolor{black}{We use the classical definition of a FLOP, counting $+$, $-$, $\times$, and $\div$, and excluding the floor and ceiling evaluations. Sorting algorithms will be treated separately.}

Consider the $y$-intercept of $\color{black} G_N \color{black}(x,\rho)$:
\begin{equation*}
\color{black} G_N \color{black}(0,\rho) =  \frac{N(N+1)}{2}\rho - \frac{N}{2}  - \sum_{k=1}^N \lfloor  k \rho \rfloor.
\end{equation*}
Computing the sum $\sum_{k=1}^N \lfloor  k \rho \rfloor$ requires $N$ multiplications and $N-1$ additions, totaling $2N-1$ operations. Adding the remaining terms contributes seven additional operations, yielding a total of $2N+6$. 

Each discontinuity occurs at $x=\lceil k \rho \rceil - k \rho $, requiring three operations per discontinuity \textcolor{black}{(assuming $k\rho$ is computed twice rather than stored)}. For $N$ discontinuities \textcolor{black}{(for irrational $\rho$, there are $N$ discontinuities, and for rational $\rho$ we track duplicates)}, this adds $3N$ operations. The endpoints of each branch are computed via $b_j = a_j + N \Delta x[j]$ (two operations) and $a_j = b_{j-1}-1$ (one operation) \textcolor{black}{ for branches $2, 3, \cdots N+1$, adding another $3N$ operations. The first branch only requires computation of $b_1=a_1+N\Delta x[1]$ (two operations) as the $y$-intercept is already computed.} \textcolor{black}{Finally, computing $\Delta x_{\texttt{list}}$ via subtraction of consecutive points requires $N+1$ operations in the irrational case --- possibly fewer in the rational case depending on its relationship with $N$. Therefore,} defining the full structure of $\color{black} G_N \color{black}(x, \rho)$ requires \textcolor{black}{at most} $\color{black} (2N+6) + (3N)+(3N+2) + (N+1) = 9N+9$ operations.

The overall order of the algorithm will depend on the choice of sorting algorithm. Our benchmarking in Section~\ref{sec:benchmark} was performed in MATLAB~
\cite{MATLAB:R2025b}; \textcolor{black}{its \texttt{sort} function uses Quicksort~\cite{McKeeman2004Sort} which has a runtime of $O(N \log N)$ on average~\cite{sedgewick1978implementing}. Therefore, the overall order of the algorithm, accounting for the initial computation of the points and the sorting is $O(N + N\log N) = O(N \log N)$. However, in our tests up to $N=10^7$, the runtime was effectively linear (see Section~\ref{sec:benchmark}).}

\subsection{Exact implementation of the $\text{pdf}$}
The DTA (Algorithm~\ref{alg:main}) provides an exact definition of the branch start and endpoints, $\left\lbrace a_i\right\rbrace$ and $\left\lbrace b_i\right\rbrace$ respectively. This collection of $2N+2$ points (may be fewer for rational $\rho$, depending on $N$) partitions the range of $\color{black} G_N \color{black}(x,\rho)$ into bins. Then we just count the branches that \textcolor{black}{overlap} each bin. This count may be accomplished through an implementation of the sweep line method which is $O(N \log(N))$~\cite{shamos1976geometric}. A naive implementation in which every interval is checked against every segment would be $O(N^2)$.

The statistical measures of $\text{pdf}(y)$, closely related to the $L^p$-norms of $\color{black} G_N \color{black}$, are typically studied in the context of their growth as a function of $N$. Thus, efficient algorithms are of utmost importance for a feasible study. 

\subsubsection*{Support}
The extrema of $\color{black} G_N \color{black}$, which define the support of $\text{pdf}$, are easily found using built-in max/min functions.

\subsubsection*{Variance}

The variance, as given in Equation~\ref{eq:variance}, can be easily simplified into an exact discrete equation by considering both the piecewise-constant structure of $\text{pdf}(y)$ and its compact support. Thus, we have

\begin{equation} \label{eq:variance_discrete}
    \sigma^2 = \int_{-\infty}^\infty y^2 \text{pdf}(y) \,dy = \sum_{j=1}^m \frac{c_j}{3} \left( y_{j+1}^3 - y_j^3\right)
\end{equation}
where $m$ is the number of piecewise-constant line segments, $c_j$ is the constant value of $\text{pdf}(y)$ on the segment, and $y_{j}$ and $y_{j+1}$ are the starting and ending $y$-values of that segment, respectively. $\sigma^2 = ||\color{black} G_N \color{black}||_2^2$ as established in Equation~\ref{eq:variance}. This computation is $O(N)$.

\subsubsection*{Kurtosis}

Kurtosis, as given in Equation~\ref{eq:kurtosis}, can be computed in a similar manner to the variance:
\begin{equation} \label{eq:kurtosis_discrete}
    \kappa = \frac{1}{\sigma^4}  \int_{-\infty}^{\infty}y^4 \text{pdf}(y)\,dy =\frac{1}{\sigma^4} \sum_{j=1}^m \frac{c_j}{5} \left( y_{j+1}^5 - y_j^5\right).
\end{equation}
Using the variance $\sigma^2$ previously computed, this computation is $O(N)$.

\section{Algorithm benchmarking} \label{sec:benchmark}
In Section~\ref{sec:numerical_implementaion}, we showed that our algorithm achieves maximum error of $O(N \varepsilon_{\text{mach}})$, and \textcolor{black}{that the number of floating point operations required to construct $G_N$ grows linearly with $N$}. Now we must show that these claims hold in practical application. All tests are run in MATLAB R2025b~
\cite{MATLAB:R2025b}.

\subsection{Time complexity} \label{sec:time}
The naive method samples $\color{black} G_N \color{black}(x,\rho)$ at some subset of points in $[0,1)$. As discussed in Section~\ref{sec:naive}, uniform sampling is unable to fully resolve the function, \textcolor{black}{but} for the purposes of making a time comparison, we take a very sparse sample: $N$ points uniform in $[0,1)$. \textcolor{black}{Of course, $N$ samples for up to $N+1$ branches are far too sparse to even roughly resolve the spatial structure of $G_N$, but good spatial resolution is computationally prohibitive for $N$ in the regime considered here.}

\begin{figure}[H]
\centering
\subfigure[Quadratic complexity for naive algorithms.]{
    \includegraphics[width=0.45\linewidth]{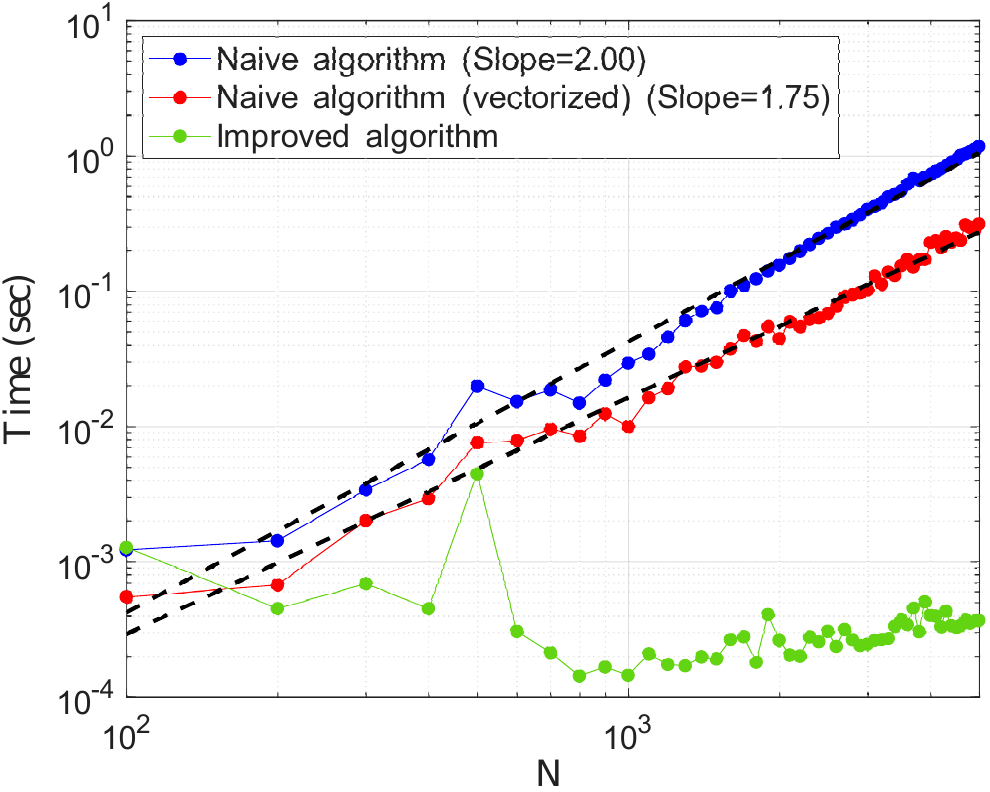}}%
\hfill
\subfigure[Linear complexity for the improved algorithm.]{
    \includegraphics[width=0.45\linewidth]{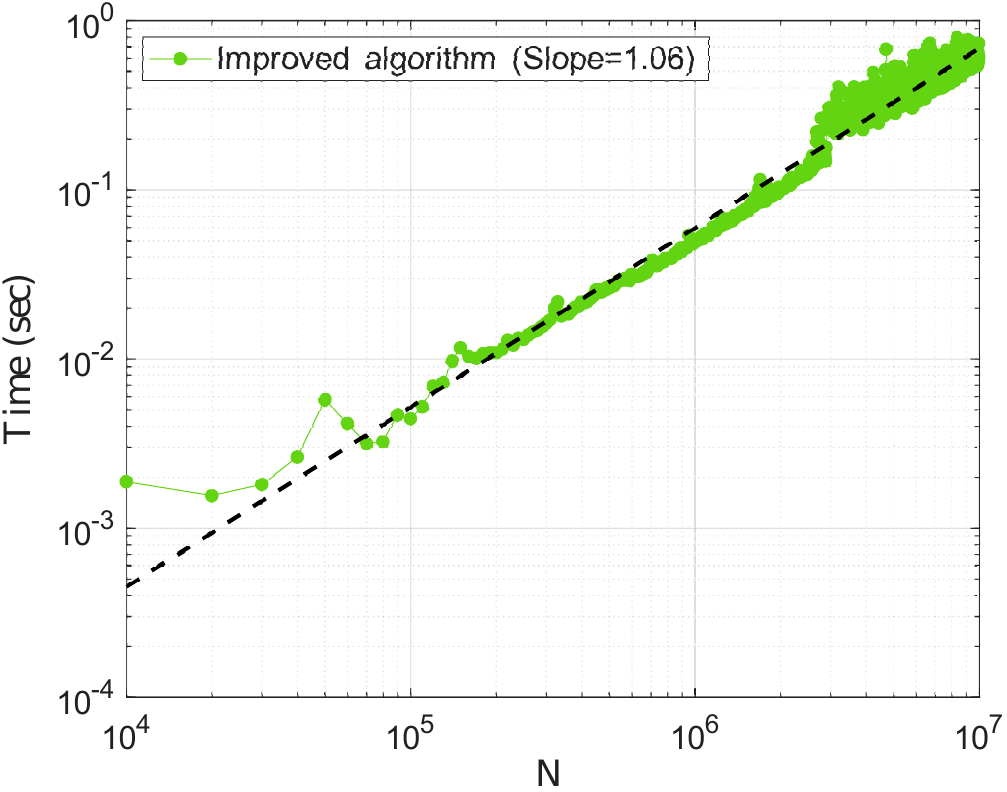} }
\caption{Comparison of time complexities for the algorithms.}
\label{fig:time_complexity_comparison}
\end{figure}

Figure~\ref{fig:time_complexity_comparison} shows runtime results for $\rho = \frac{\sqrt{5}-1}{2}$ for the naive algorithm (Algorithm~\ref{alg:naive}), its vectorized version, and the DTA (Algorithm~\ref{alg:main}). The left panel confirms that the naive algorithm is $O(N^2)$, improved to $O(N^{1.75})$ through vectorization. The right panel shows that the DTA runs in $O(N)$ time (represented by a linear fit with slope of 1 and $R^2=0.97$). Similar results hold across choices of $\rho$; the same analysis was run with $\rho =L$ (slope of 1.07, $R^2=0.97$).

\subsection{Error analysis} \label{sec:error}
\begin{figure}[H]
\centering
    \includegraphics[width=0.48\linewidth]{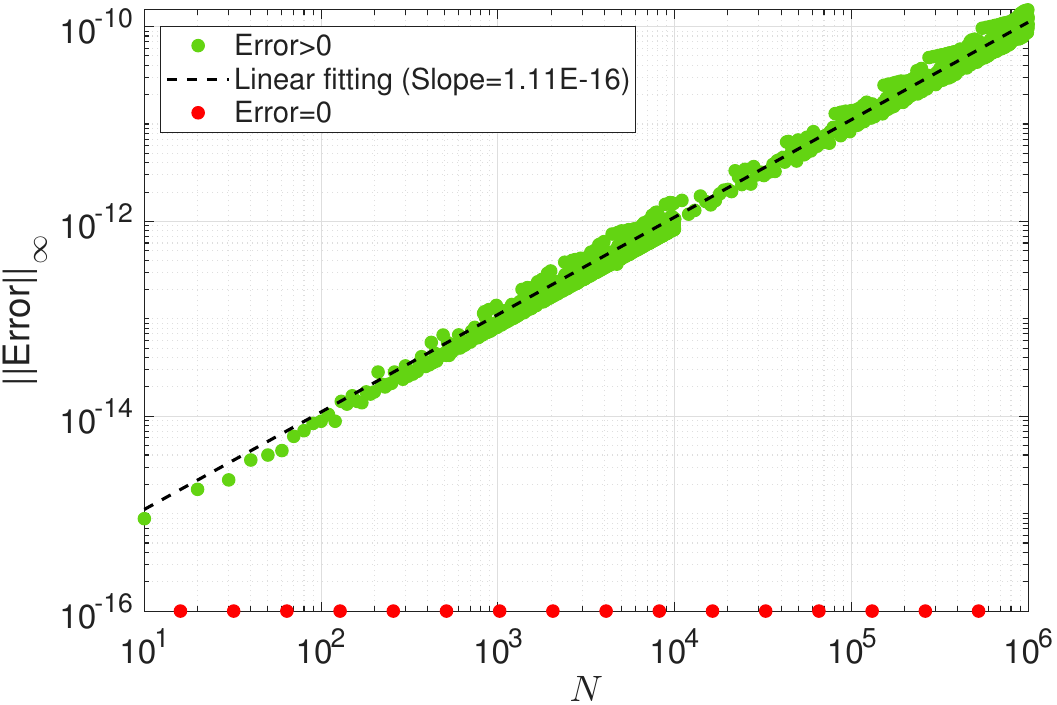} 
\caption{The maximum error ($\infty$-norm) of the algorithm grows with the number of terms as $\frac{1}{2}N \varepsilon_{\text{mach}}$.}
\label{fig:accuracy}
\end{figure}

\textcolor{black}{When $\rho$ is rational}, $\color{black} G_N \color{black}$ exhibits a predictable structure for particular $N$. When $\rho = \frac{p}{q}$ \textcolor{black}{with $\gcd(p,q)=1$} and $N = kq$, $\color{black} G_N \color{black}$ forms $q$ branches spanning $ [-\frac{k}{2}, \frac{k}{2} ]$ and $x \in [\frac{r-1}{q}, \frac{r}{q}]$ for $r = 1, 2, \dots, q$ (Changing $p$ does not change the behavior for $N=kq$). Vertical jumps of height $k$ occur between branches. To evaluate the accuracy of the algorithm, we set $\rho = \frac{1}{q}$ and $N = q$. The infinity norm of the error is the maximal endpoint deviation from $\pm \frac{1}{2}$. Results in Figure~\ref{fig:accuracy} show that errors in $\color{black} G_N \color{black}$ grow proportionally to $\frac{1}{2}N \varepsilon_{\text{mach}}$, as confirmed by a linear fit with $R^2 = 0.95$.  Errors are exactly zero for $q = 2^j, j \in \mathbb{N} $, corresponding to $\rho$ values exactly representable in double-precision. Of course these zero errors can only be a feature of rational $\rho$.

\section{Discussion}
\label{sec:discussion}
The improved algorithm overcomes the $O(N^2)$ computational barrier, and the numerical stability barriers, bringing the study of the discrepancy of irrational rotations firmly into the regime of experimental mathematics. We now present some examples of what can be studied using this tool. We focus on the discrepancy, \textcolor{black}{i.e.,} $\operatorname{Supp}(\text{pdf})$ (twice $||\color{black} G_N \color{black}||_{\infty}$), as it is \textcolor{black}{a} common metric in the theoretical literature. The variance (equivalently, $||\color{black} G_N \color{black}||_{2}^2$) and kurtosis exhibit behavior highly correlated with $||\color{black} G_N \color{black}||_{\infty}$, as all three metrics respond to the branch realignments in related ways, e.g. \textcolor{black}{appear to} have local extrema at \textcolor{black}{convergent denominators} and exhibit self-similarity with respect to \textcolor{black}{their} multiples (see Appendix~\ref{apx:stats_plots}).

\subsection{Patterns in the \textcolor{black}{discrepancy}}
It is a foundational result in number theory that the behavior of $\color{black} \|G_N (\cdot,\rho)\|_{\infty}$ \textcolor{black}{(half the discrepancy)} is governed by the partial quotients of $\rho$; in particular, \textcolor{black}{discrepancy is low} when $N=q_n$ is a \textcolor{black}{convergent denominator} of $\rho$~\cite{kuipers1974uniform, baxa1994minimum}. Therefore, in order to capture the structure of \textcolor{black}{the discrepancy} as $N$ increases without sampling too many $N$, we use the \textcolor{black}{convergent denominators} and their multiples as anchors and sample logarithmically between them. In order to find the local maxima numerically, a different sampling strategy would be required.

\begin{figure}[H]
\centering
\subfigure[$\rho =L= \color{black} \lbrack 0; 9, 11, 99, 1, 10, 9, 10^{12}-1 , 1, \dots \rbrack$]{
    \includegraphics[width=0.4\linewidth]{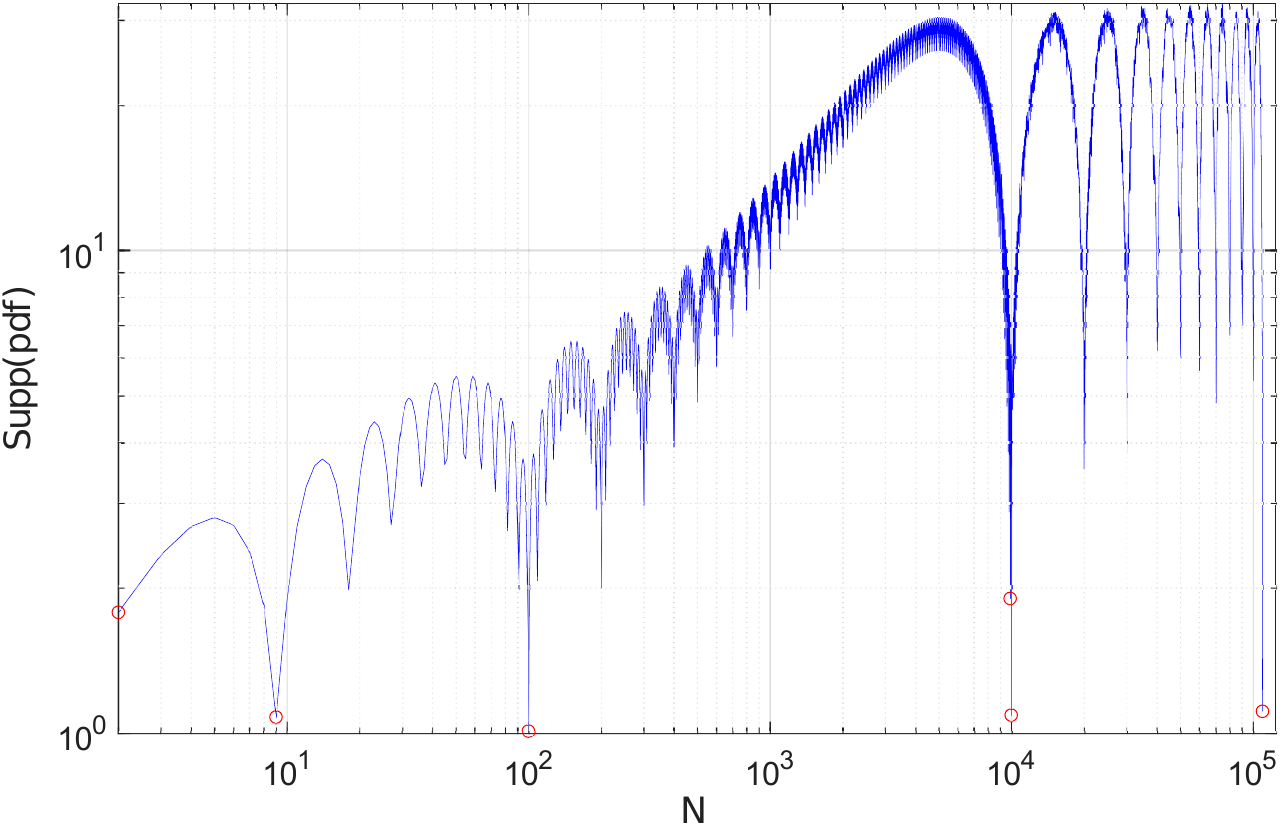}}%
\subfigure[$\rho=\pi-3= \color{black} \lbrack 0; 7, 15, 1, 292, 1, 1, 1, 2, \dots \rbrack$]{
    \includegraphics[width=0.4\linewidth]{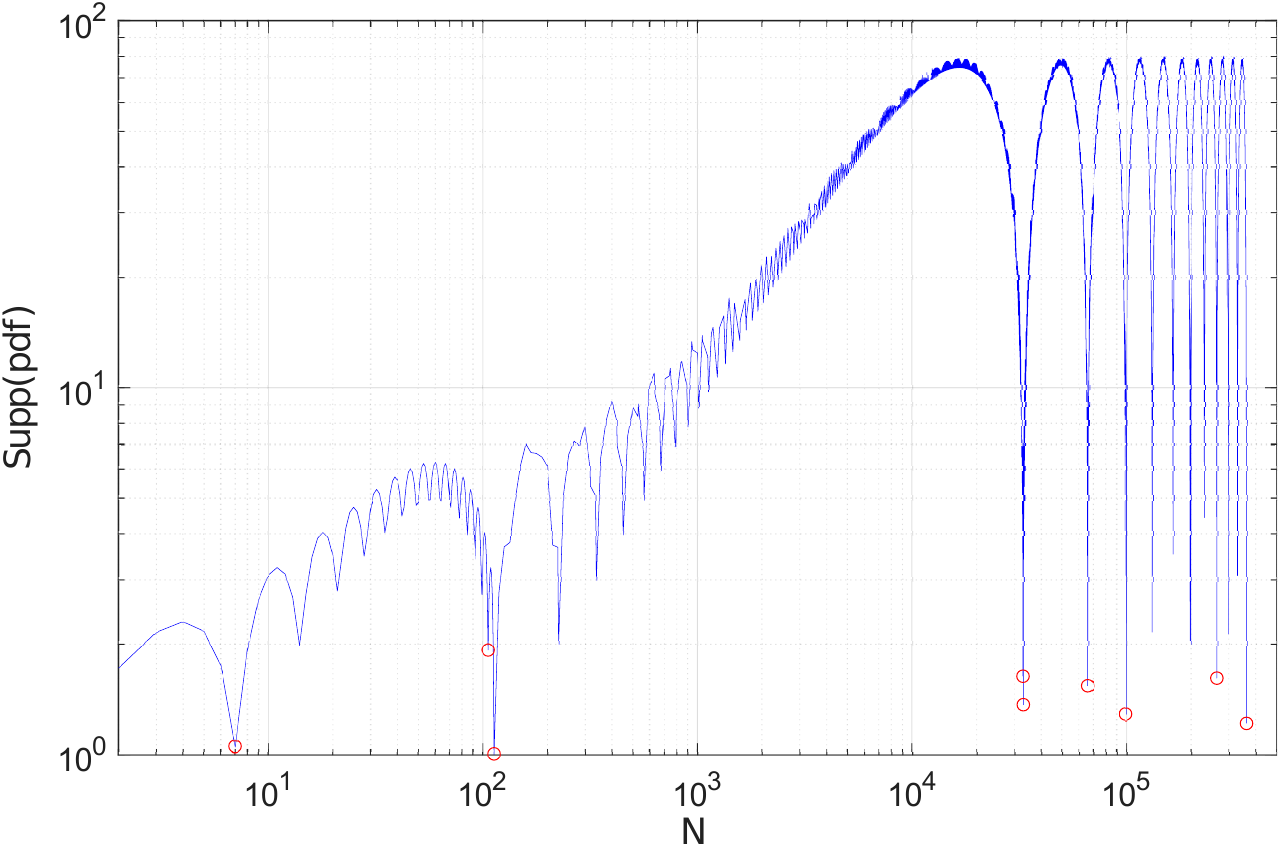}} 
    \subfigure[$\rho=e-2= \color{black} \lbrack 0; 1, 2, 1, 1, 4, 1, 1, 6, \dots \rbrack$]{
    \includegraphics[width=0.4\linewidth]{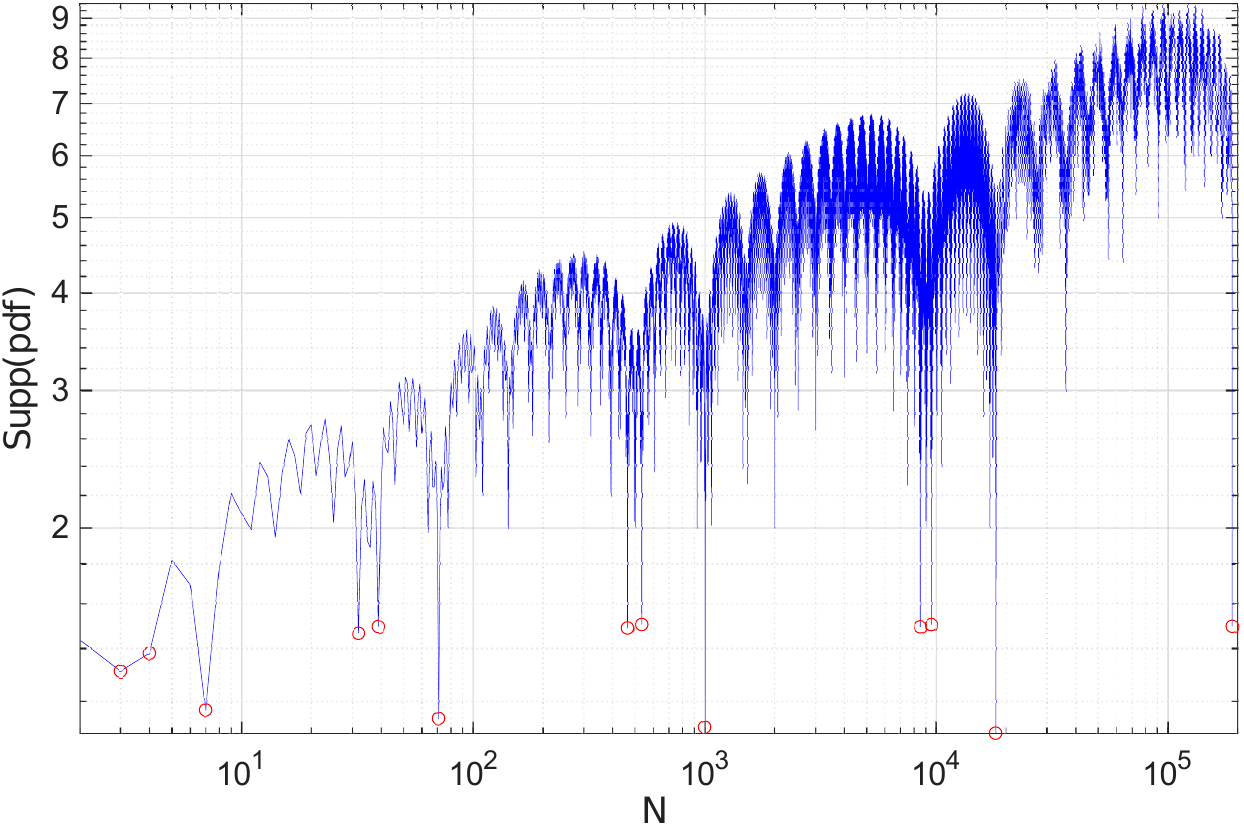}}%
\subfigure[$\rho=\frac{\pi^2}{6}-1=\color{black} \lbrack 0; 1, 1, 1, 4, 2, 4, 7, 1, 4, \dots \rbrack$]{
    \includegraphics[width=0.4\linewidth]{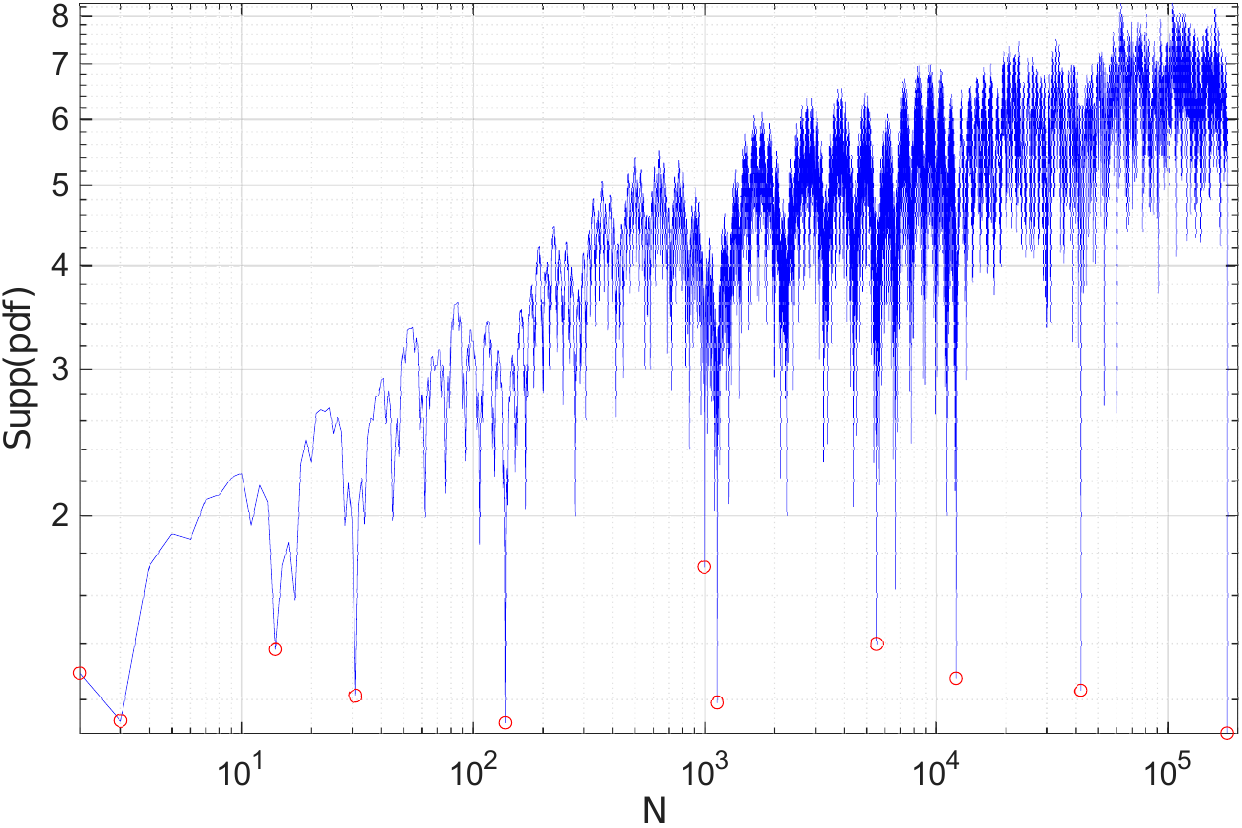}}
        \subfigure[$\rho=\sqrt{2}-1 \color{black}= \lbrack  0; \overline{2} \rbrack$]{
    \includegraphics[width=0.4\linewidth]{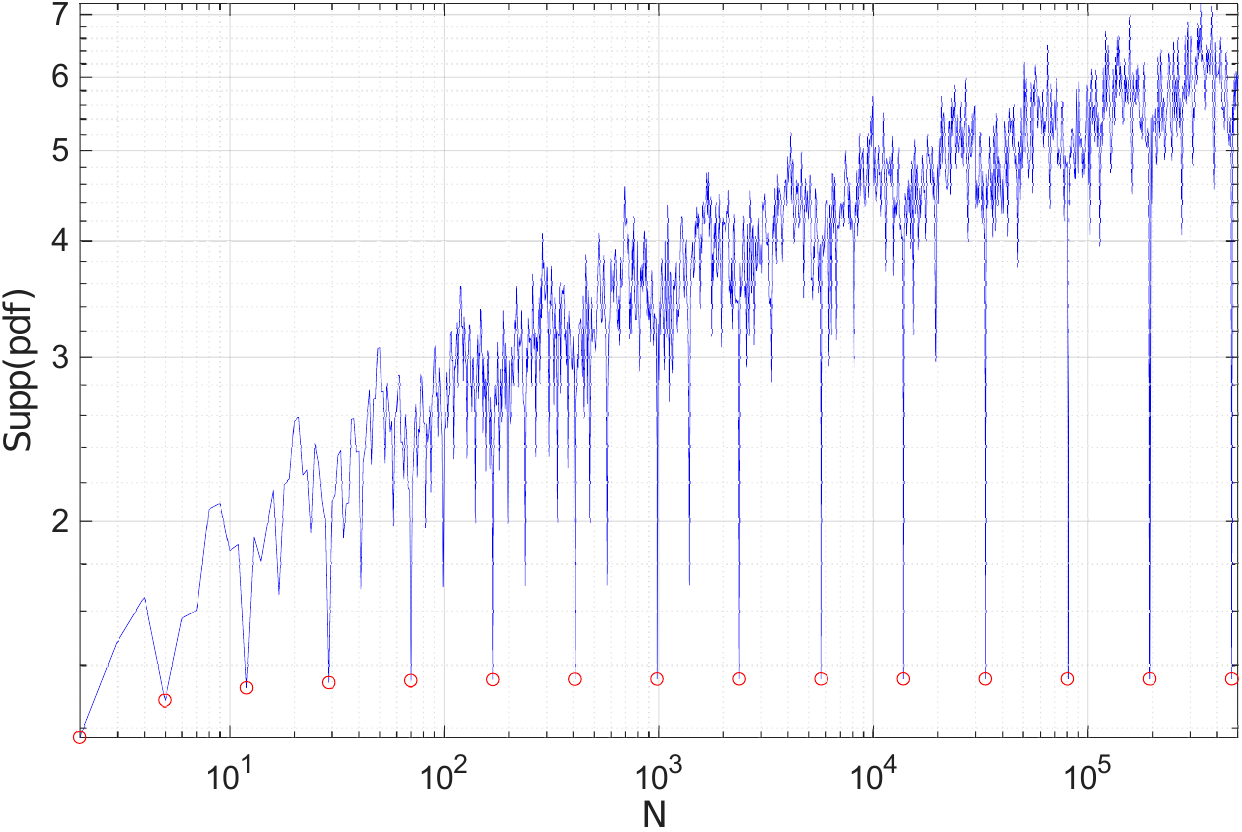}}%
   \subfigure[$\rho=\alpha_{ST} \color{black}= \lbrack  0; 2, 40, 40,  \overline{2} \rbrack$]{
    \includegraphics[width=0.4\linewidth]{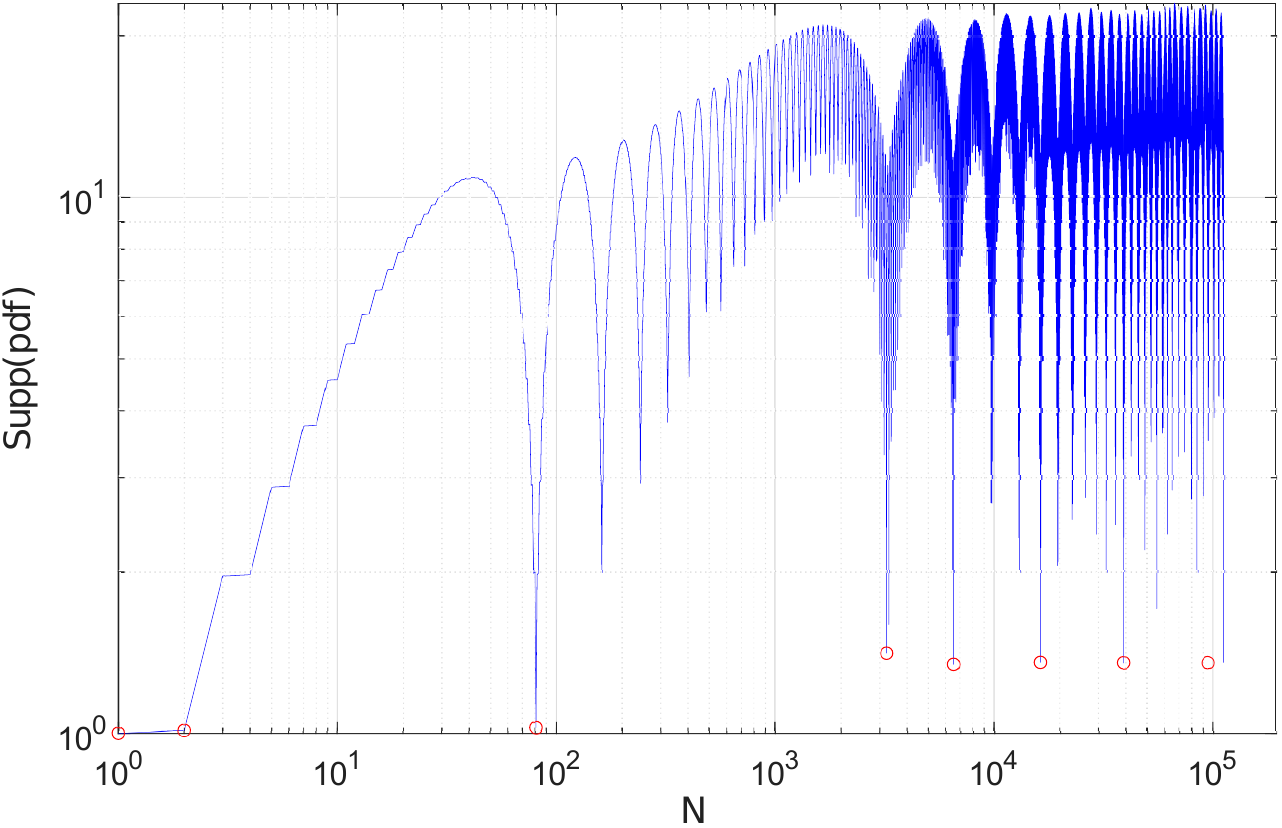}} 
\caption{Plots of $\operatorname{Supp}(\text{pdf}(y))$ for various classes of irrational numbers for $N=2, 3, \dots$, with the \textcolor{black}{convergent denominators} plotted in red. (Top row) Transcendental numbers with a large partial quotient. (Second row) Transcendental numbers with a medium-sized partial quotient. (Third row) Quadratic irrationals: on the left is $\sqrt{2}\color{black}-1\color{black}=\lbrack \color{black} 0;  \overline{2} \rbrack$, and on the right is $\alpha_{ST}= \lbrack \color{black} 0; \color{black}2, 40, 40, \color{black} \overline{2} \rbrack$ studied by Shimaru \& Takashima~\cite{shimaru2018discrepancies}.}
\label{fig:support_diophantine}
\end{figure}

Figure~\ref{fig:support_diophantine} shows $\operatorname{Supp}(\text{pdf})$ on a $\log$-$\log$ plot for a few classes of irrationals, with the \textcolor{black}{convergent denominators} plotted as red dots. Various types of self-similarity can be observed. The top row shows $\rho=L = \color{black} \lbrack 0; 9, 11, 99, 1, 10, 9, 999999999999, 1, \dots \rbrack$ (Liouville's constant) and $\rho=\pi\color{black} -3\color{black}= \color{black} \lbrack 0; 7, 15, 1, 292, 1, 1, 1, 2, \dots \rbrack$, whose support functions in this regime are characterized by very large hills preceding the large \textcolor{black}{convergent denominator} (9,909 for $L$ and 33,102 for $\pi\color{black} -3\color{black}$), as observed by Vinson~\cite{vinson2001partial} for $\rho=\pi$ and by Shimaru \& Takashima~\cite{shimaru2018discrepancies} and Setokuchi \& Takashima~\cite{setokuchi2014discrepancies} for quadratic irrationals with similar properties. The second row shows two transcendentals without particularly high partial quotients \textcolor{black}{in this regime}, $\rho=e-2 = \color{black} \lbrack 0; 1, 2, 1, 1, 4, 1, 1, 6, \dots \rbrack$ and $\rho=\frac{\pi^2}{6}-1 = \color{black} \lbrack 0; 1, 1, 1, 4, 2, 4, 7, 1, 4, \dots \rbrack$. The last row shows two quadratic irrationals, which are characterized by eventually-repeating \textcolor{black}{partial quotients}. $\rho=\sqrt{2} \color{black} -1\color{black} = \lbrack \color{black} 0; \overline{2} \color{black}\rbrack$ cannot sustain large peaks due to the regularity of its partial quotients. $\rho=\alpha_{ST}= \lbrack \color{black} 0;\color{black}2, 40, 40, \overline{2} \rbrack$, a quadratic irrational studied by Shimaru \& Takashima~\cite{shimaru2018discrepancies}, has some large early hills then settles into a more regular structure without large peaks. \textcolor{black}{The discrepancy seems to exhibit local minima when $N=q_n$}. 

\subsubsection*{\color{black} Self-similarity}

\begin{figure}[H]
\centering
        \subfigure[$\rho=\pi-3$]{
    \includegraphics[width=0.75\linewidth]{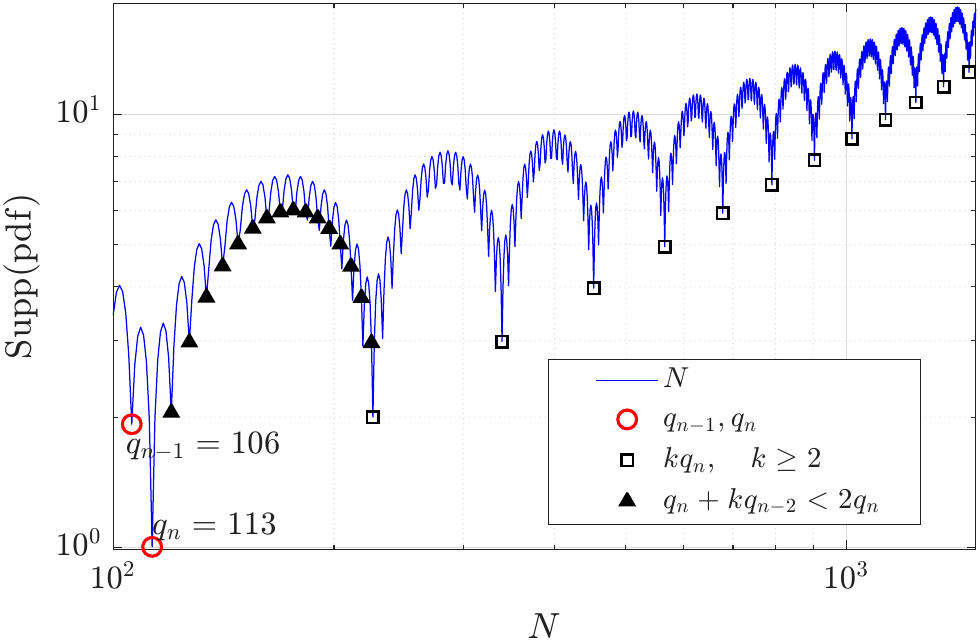}} \\
      \subfigure[$\rho=\alpha_{ST}$]{
    \includegraphics[width=0.75\linewidth]{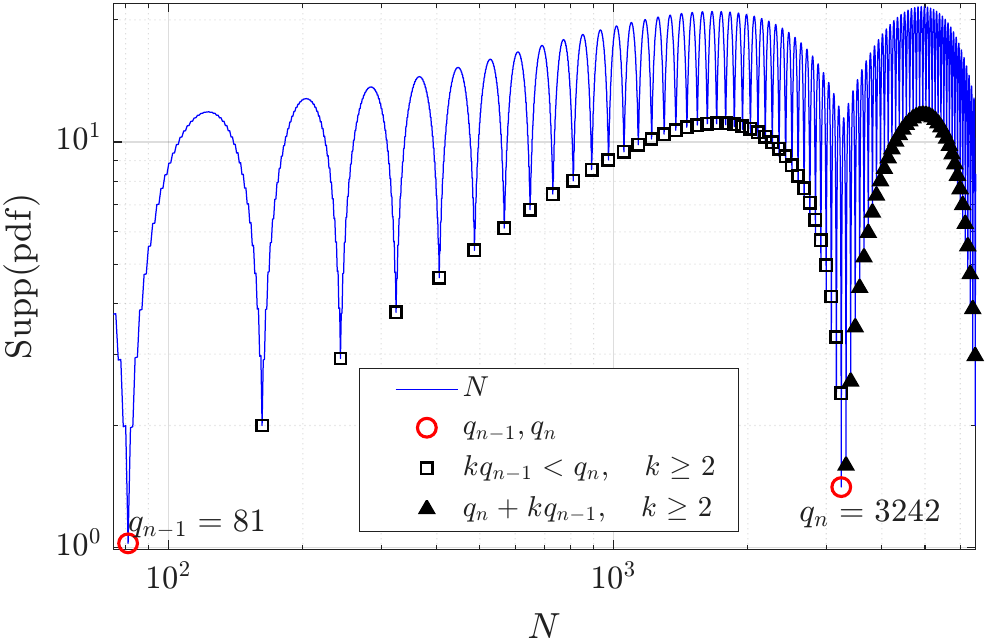}}%
\caption{The figure shows zoomed-in versions of some plots in Figure~\ref{fig:support_diophantine} with some special values of $N$ marked.}
\label{fig:selfsim}
\end{figure}

\textcolor{black}{Figure~\ref{fig:selfsim} shows the details of some of what we call ``self-similar'' behavior for two irrationals: $\rho=\pi-3$, whose partial quotients do not repeat with any known pattern, and $\rho=\alpha_{ST}= \lbrack 0; 2, 40, 40, \overline{2} \rbrack$, a quadratic irrational (partial quotients are eventually periodic, and in this case, eventually constant). In both cases, the plots show nested arches organized around convergent denominators and their multiples. These structures are not exactly periodic, however, as the locations of their local minima seem to shift as $N$ increases. For example, in the top panel  ($\rho=\pi-3$), the local minima initially occur at $N=113+7k = q_3+q_1k$, but after several repetitions, they shift by one unit.}

\textcolor{black}{We also note an exact expression for some of these points of interest. When $N=kq_n$ for $1 \leq k \leq a_{n+1}$, where $a_{n+1}$ is the $n+1$st partial quotient~\cite{baxa1994minimum},
\begin{equation}
D_{kq_n}(\rho)=k\left[1-\left((k-2)q_n+1\right) \left|p_n-q_n\rho\right|\right].
\end{equation}
When $k=1$, the values are the convergent denominators of $\rho$, i.e. the red circles denoted in Figures~\ref{fig:support_diophantine} and~\ref{fig:selfsim}.}

\subsection{Recurring patterns in the $\text{pdf}$}

\subsubsection*{Trapezoids}

As discussed in the previous section,  $||\color{black} G_N \color{black}(\cdot,\rho)||_{\infty}$ (half of $\text{Supp}(\text{pdf})$) is \textcolor{black}{small} when $N=q_n$ is a \textcolor{black}{convergent denominator}~\cite{kuipers1974uniform}. These are the points where $\rho$ is best approximated by a rational number $\frac{p_n}{q_n}$. This manifests in the plot of $\color{black} G_N \color{black}$ as branches that are highly aligned. An earlier version of this algorithm developed by the present author was used by Ralston et al. to provide supporting numerical evidence for their ``trapezoid theorem,'' which states that this alignment causes the $\text{pdf}$ to reach its maximum possible height of 1 \textcolor{black}{when} $N=q_n$~\cite{veerman_BS}. 

Here we provide some brief geometric intuition for the theorem. For rational $\rho = \frac{p}{q}$ \textcolor{black}{with $\gcd(p,q)=1$ and} $N=q$, \textcolor{black}{the function} $\color{black} G_N \color{black}$ forms $q\color{black}$ perfectly aligned branches spanning $\left[-\frac{1}{2}, \frac{1}{2}\right]$. The corresponding $\text{pdf}$ is the rectangular uniform distribution $\left[-\frac{1}{2}, \frac{1}{2}\right] \times [0,1]$. When $\rho$ is irrational and $N=q_n$ is a convergent denominator, $\color{black} G_N \color{black}$ is characterized by $\color{black} q_n+1$ branches, \textcolor{black}{$q_n$} of which overlap on a smaller subset of $\left[-\frac{1}{2}, \frac{1}{2}\right]$, with some shifted a bit on either end. The height still reaches 1, but the length of the top is now shorter than 1. This causes the base to widen, as the $\text{pdf}$ is a distribution and its area must always equal 1. The result is a trapezoid whose dimensions depend on the quality of the rational approximation $\frac{p_n}{q_n}$ to $\rho$. The shape is not a ``true'' trapezoid in that the sides are made up of small piecewise constant line segments.

\subsubsection*{Trapezoids with spikes}
Numerical experiments show particular patterns for $\rho$ when it is well-approximated by rational $\frac{p_n}{q_n}$ and $N=kq_n$ for $k=1, 2, \dots m$. In this case, we see the distribution develop $k-1$ spikes, starting from the center. The \textcolor{black}{density} between the spikes remains completely flat until a particular $m$ when the shape begins to degenerate into oscillatory noise. How quickly this degeneration happens \textcolor{black}{appears to} depend on how well-approximated $\rho$ is by $\frac{p_n}{q_n}$. We can quantify the quality of the approximation by the local irrationality exponent:
\begin{equation}
\left| \rho - \frac{p_n}{q_n}\right| \color{black} \approx \color{black} \frac{1}{q_n^{\mu_{\color{black} n}}} \implies \mu_n \color{black}\approx\color{black} -\frac{\log\left| \rho - \frac{p_n}{q_n} \right|}{\log(q_n)}
\end{equation}
with $\mu_n \geq2$ always when $q_n$ is a \textcolor{black}{convergent denominator~\cite{Khinchin_CF}. For a Liouville irrational, $\mu_n$ is unbounded as $n$ increases. For general irrational $\rho$, individual values of $\mu_n$ may vary substantially with $n$, depending on the partial quotients.}

\begin{figure}[H]
\centering
\subfigure[$\rho=\pi-3$]{
    \includegraphics[width=0.9\linewidth]{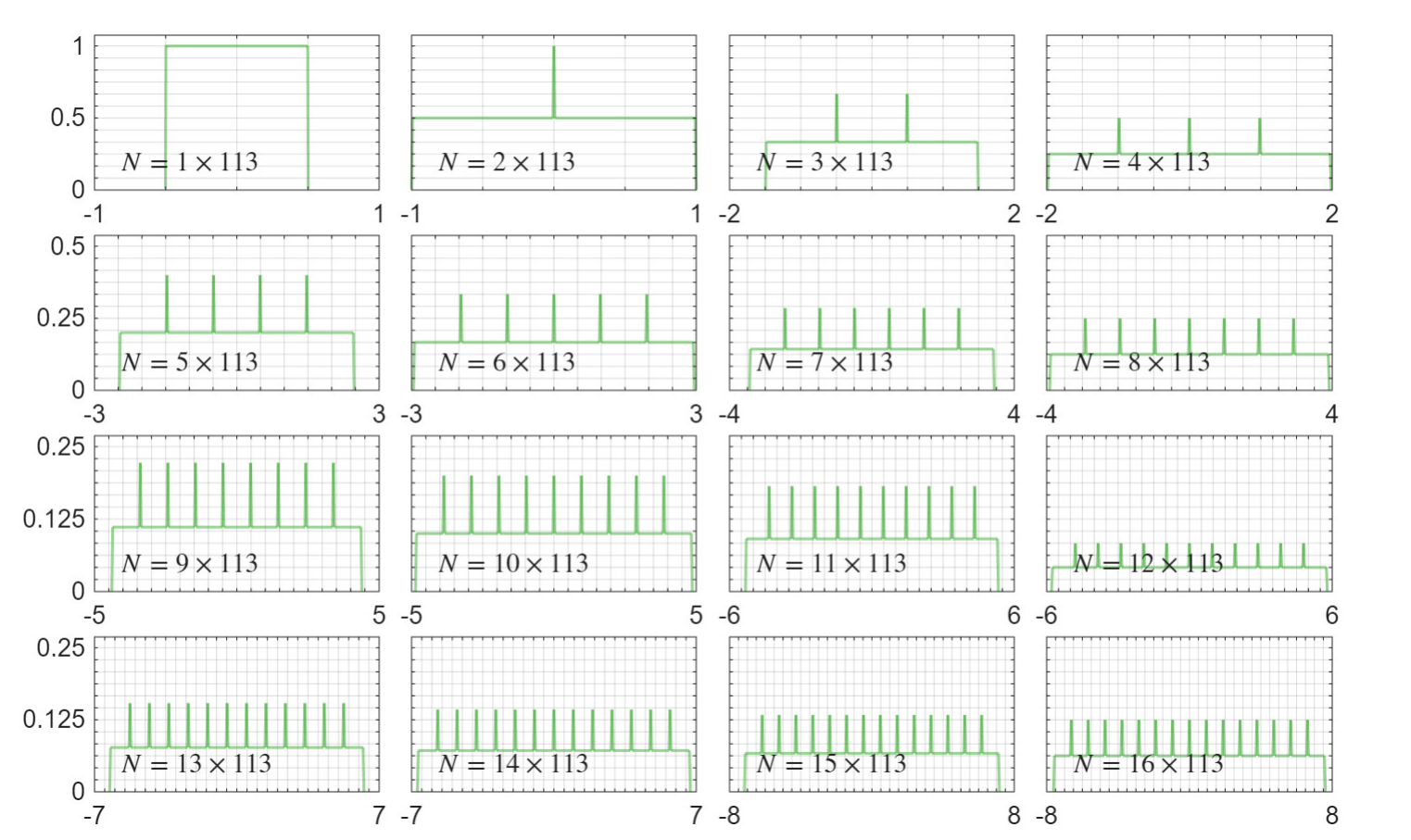}}
    \subfigure[Shimaru \& Takashima's $\alpha_{ST}= \lbrack \color{black} 0; \color{black} 2, 40, 40, \color{black} \overline{2} \color{black}\rbrack$]{
    \includegraphics[width=0.9\linewidth]{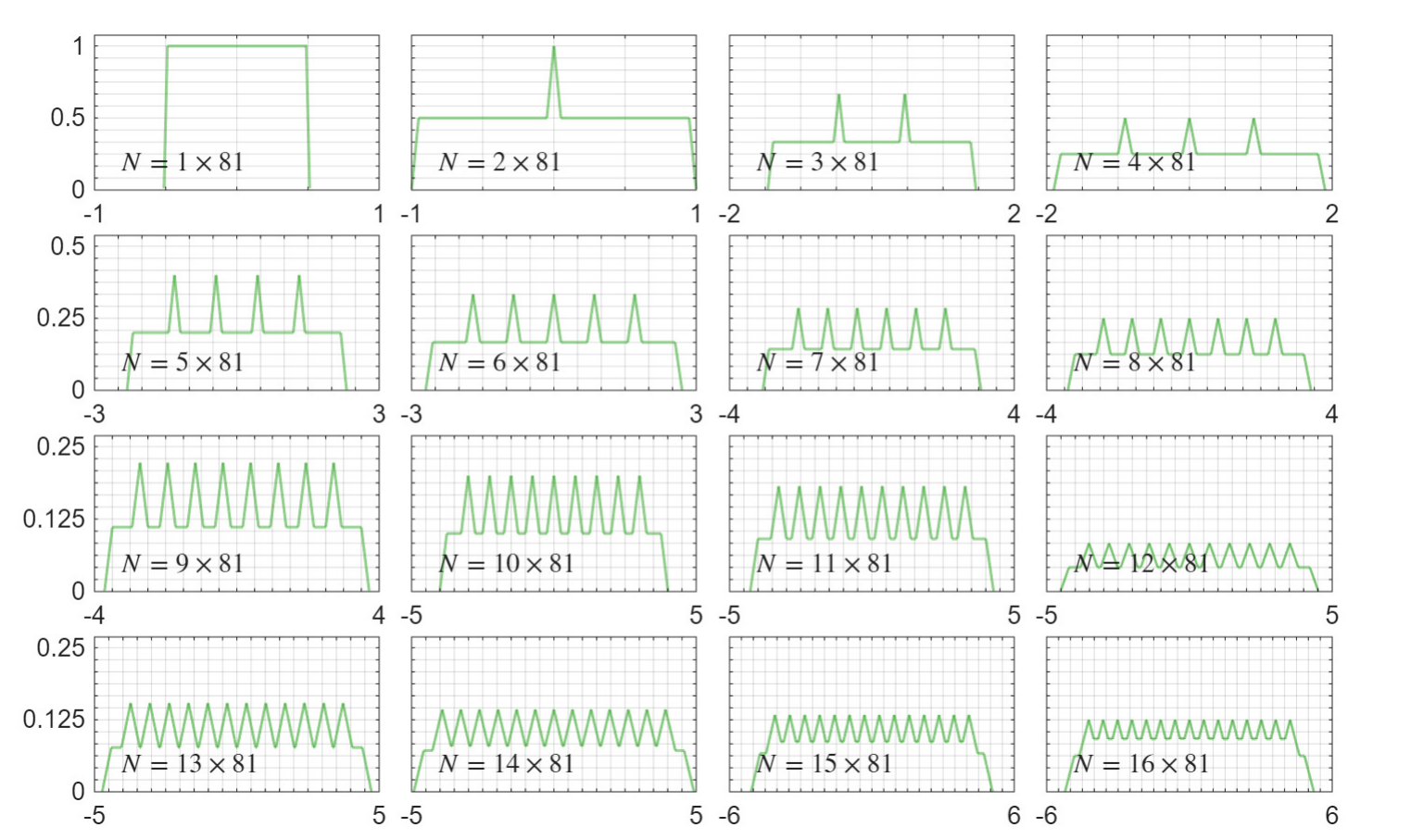}}
\caption{Multiples of \textcolor{black}{convergent denominators} for two well-approximated irrationals show trapezoids with spike patterns that eventually degrade.}
\label{fig:spikes_nondegen}
\end{figure}

Large $\mu_n$ (good rational approximations) correspond to trapezoids that are nearly rectangular, with sharply sloped sides. Such shapes are able to sustain spike formation for many multiples of $q_n$. The stability of this pattern seems to correspond with how well $\rho$ is approximated by $\frac{p_n}{q_n}$. Figure~\ref{fig:spikes_nondegen} shows the pattern for some good approximations. The first, $\left(\rho, \frac{p_n}{q_n}, \mu_n \right) = \left(\pi-3 , \frac{16}{113}, 3.2 \right)$, has a pattern that persists for many multiples of $q_n=113$. The second is $ \left(\alpha_{ST}, \frac{40}{81}, 2.8 \right)$. We see the pattern break down at $N=81 \times 14$ - the bottom of the spikes no longer reaches the top of the trapezoid.

\begin{figure}[H]
\centering
\subfigure[$\rho = e-2$]{
    \includegraphics[width=0.9\linewidth]{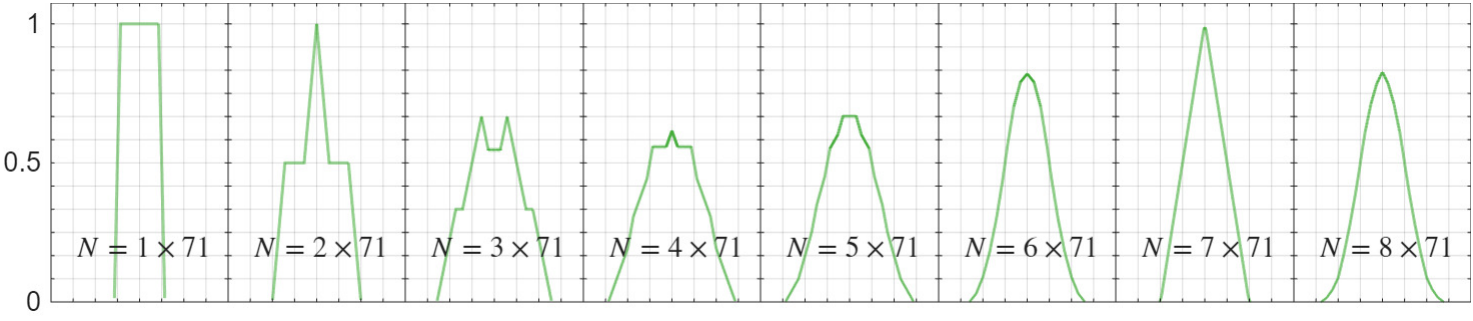}}
\subfigure[$\rho=\frac{\pi^2}{6}-1$]{
    \includegraphics[width=0.9\linewidth]{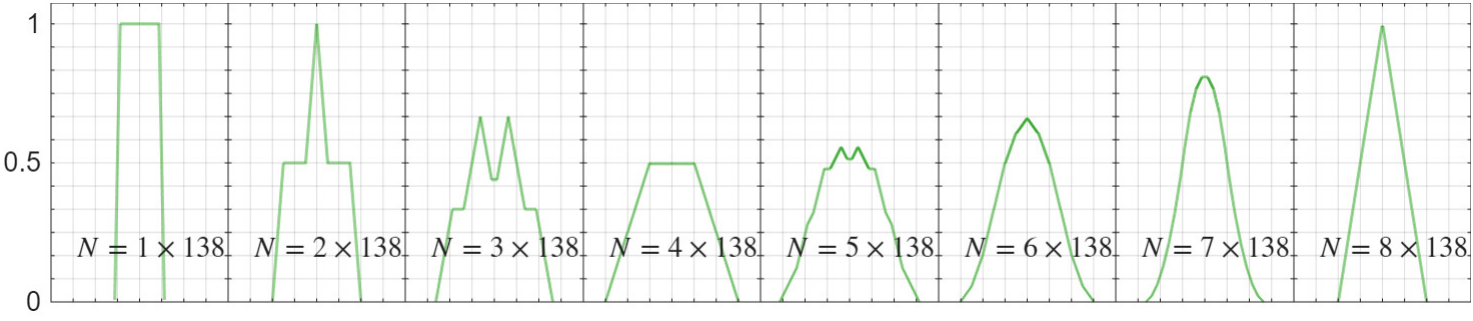}}
    \subfigure[Quadratic irrational $\rho=\lbrack \color{black} 0; \color{black}2, 2, 2, 2, 10, \color{black} \overline{2} \color{black} \rbrack$]{
    \includegraphics[width=0.9\linewidth]{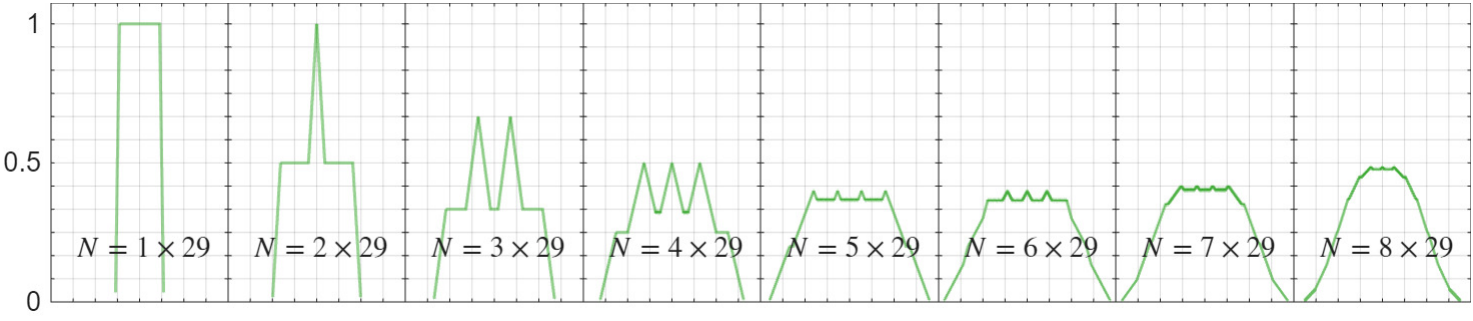}}
        \subfigure[Quadratic irrational $\rho=\lbrack \color{black} 0; \color{black}2, 2, 2, 2, 15, \color{black} \overline{2} \color{black} \rbrack$]{
    \includegraphics[width=0.9\linewidth]{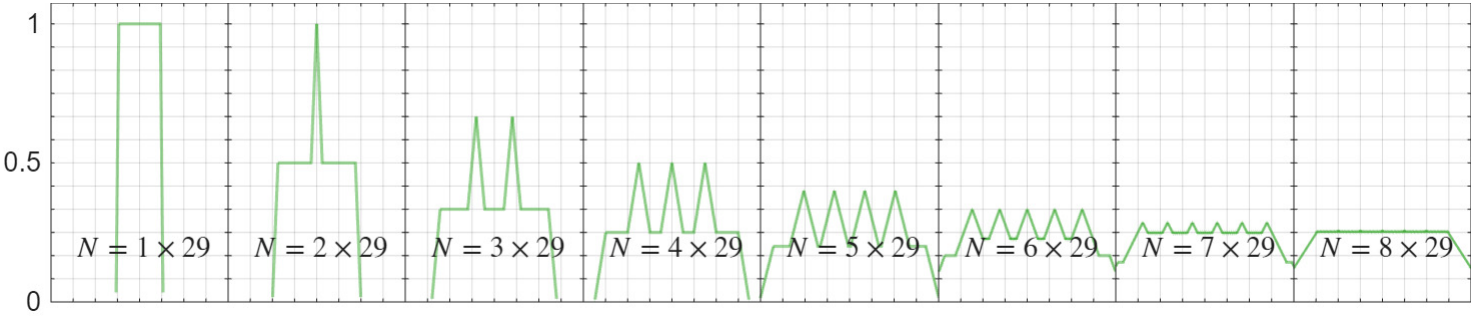}}
\caption{Multiples of \textcolor{black}{convergent denominators} for some irrationals with moderately good rational approximations show the spike structure start to degenerate early, at $k=3$ for $\rho =e-2$ and $\rho = \frac{\pi^2}{6}-1$, $k=4$ for $\rho=\lbrack \color{black} 0; \color{black}2, 2, 2, 2, 10, \color{black} \overline{2} \color{black} \rbrack$, and  $k=5$ for $\rho=\lbrack \color{black} 0; \color{black}2, 2, 2, 2, 15, \color{black} \overline{2} \color{black} \rbrack$. \textcolor{black}{All distributions are shown on a $[-2,2] \times [0,1]$ grid, with the longer tails that fall outside the grid not shown.}}
\label{fig:spikes_degen}
\end{figure}

Figure~\ref{fig:spikes_degen} shows the pattern for a few irrationals with moderately good rational approximations. We showcase the following combinations of $\left(\rho, \frac{p_n}{q_n}, \mu_n \right)$: the top two rows feature two transcendentals: (top row) $\left(e-2, \frac{51}{71}, 2.46 \right)$, (second row) $\left(\frac{\pi^2}{6}-1, \frac{89}{138}, 2.42 \right)$, while the bottom two rows showcase some quadratic irrationals $\left(\lbrack \color{black} 0; \color{black} 2, 2, 2, 2, 10,\color{black} \overline{2} \color{black}\rbrack,  \frac{12}{29}, 2.71 \right)$ and $\left(\lbrack \color{black} 0; \color{black}2, 2, 2, 2, 15, \color{black} \overline{2} \color{black}  \rbrack,  \frac{12}{29}, 2.82 \right)$. For the first two, the pattern begins to degenerate when $k=3$, the third when $k=4$, and the fourth when $k=5$.

\begin{figure}[H]
\centering
\subfigure[$\rho = \pi-3$]{
    \includegraphics[width=0.9\linewidth]{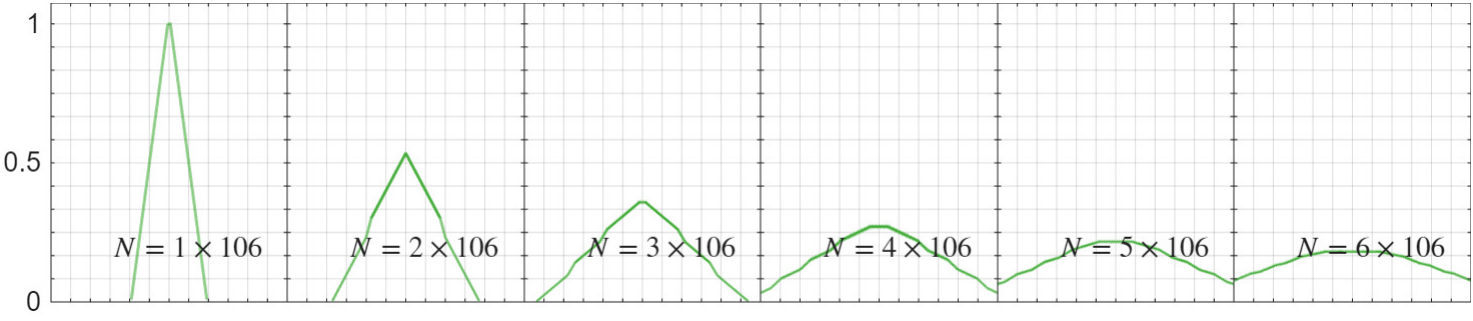}}
    \subfigure[$\rho = \pi-3$]{
    \includegraphics[width=0.9\linewidth]{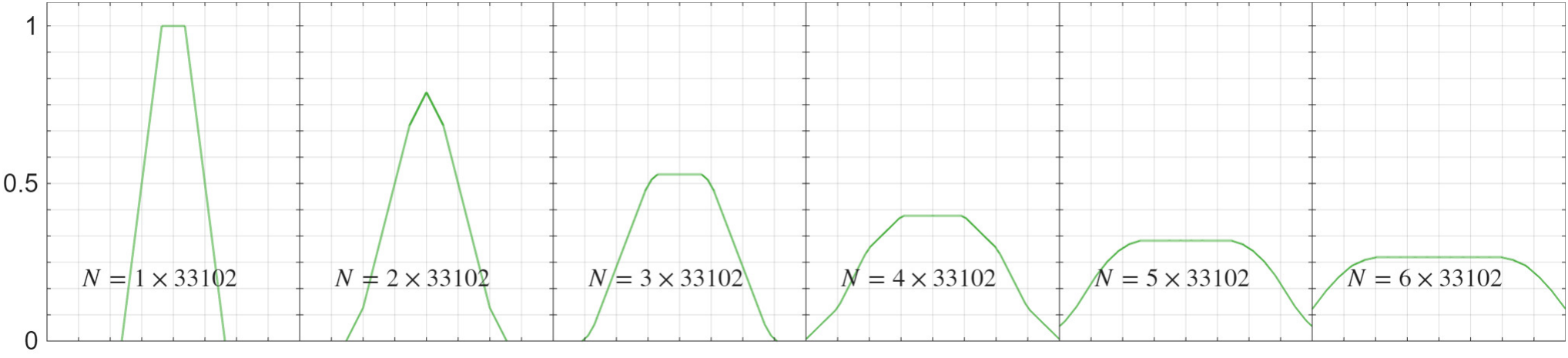}}
\subfigure[$\rho=e-2$]{
    \includegraphics[width=0.9\linewidth]{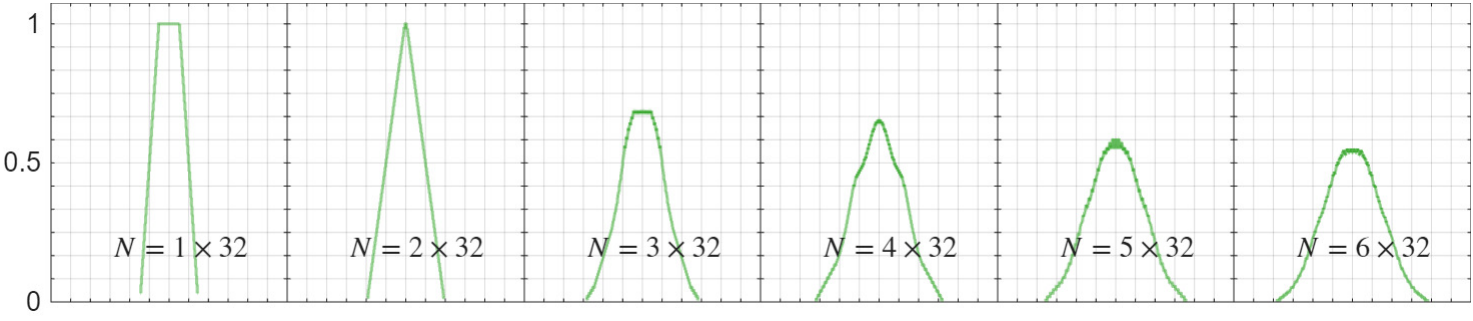}}
        \subfigure[$\rho=\sqrt{2}-1$]{
    \includegraphics[width=0.9\linewidth]{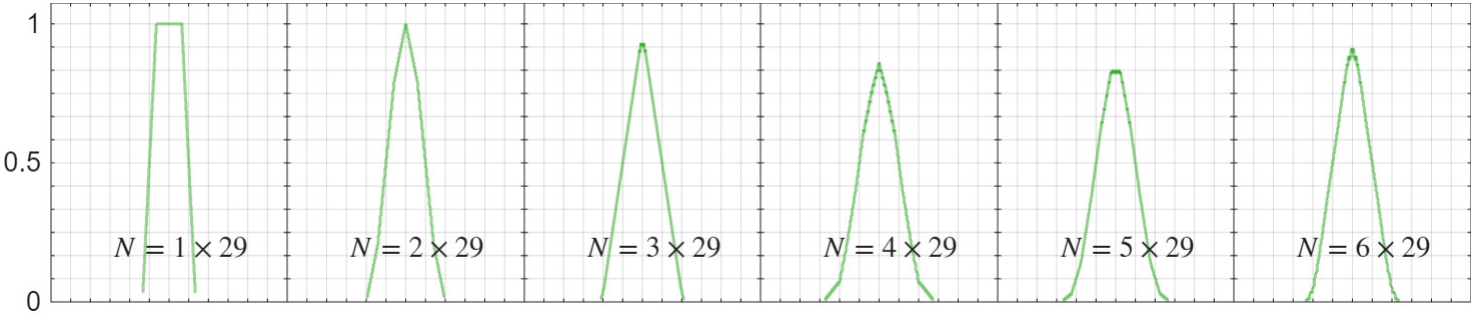}}
            \subfigure[$\rho=\sqrt{2}-1$]{
    \includegraphics[width=0.9\linewidth]{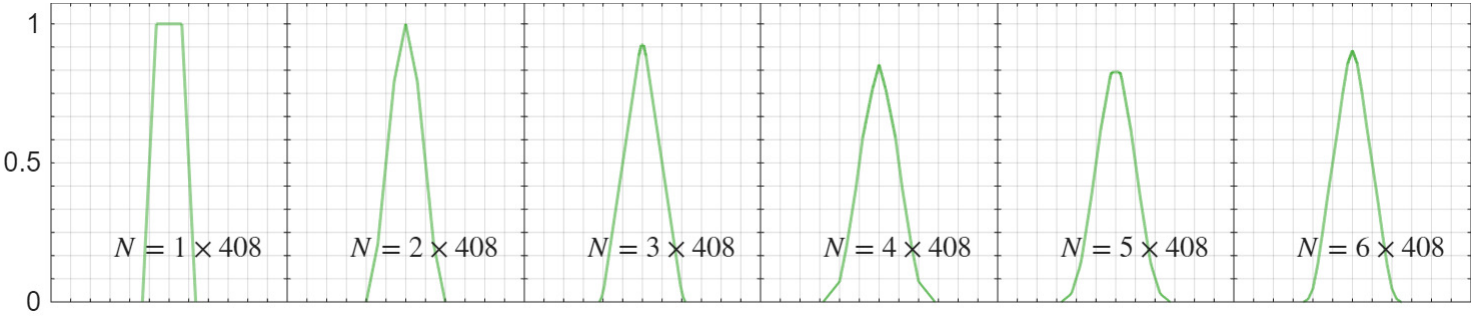}}
\caption{Multiples of \textcolor{black}{convergent denominators} for some irrationals with bad rational approximations show varied behavior.  \textcolor{black}{All distributions are shown on a $[-2,2] \times [0,1]$ grid, with the longer tails that fall outside the grid not shown.}}
\label{fig:no_spikes}
\end{figure}

When $\mu_{\color{black} n}$ is small (bad rational approximations), the trapezoid is \textcolor{black}{more pyramidal with a very short top and} unable to support the pattern for $k>1$. Some examples are shown in Figure~\ref{fig:no_spikes}.

\section{Conclusion}
We have presented an algorithm for the computational analysis of the \textcolor{black}{centered sum of fractional parts} $\color{black} G_N \color{black}(x,\rho)$. This algorithm allows rapid generation of figures \textcolor{black}{with minimal error on the order of $O(N \varepsilon_{\text{mach}})$}. By defining the structure of the sum through its discontinuities, we are able to bypass both numerical instabilities and computational bottlenecks. We showed that by using a novel discontinuity-tracking algorithm, we can fully characterize the structure of $\color{black} G_N \color{black}$ in \textcolor{black}{approximately} $O(N)$ time \textcolor{black}{for practical $N$}. This enables the rapid computation of its corresponding probability density function $\text{pdf}$ and its associated statistical properties (support, variance, and kurtosis) in $O (N \log N)$ time, which allows rapid experimental study of \textcolor{black}{the full spatial structure of $\color{black} G_N \color{black}$ as well as the discrepancy of the irrational rotation}.

We presented two examples of analyses made possible by the algorithm. First, exact plots of the length of the support, $\operatorname{Supp}(\text{pdf})$, as a function of $N$ reveal intricate self-similar patterns related to the \textcolor{black}{subsequences} of convergent \textcolor{black}{denominators}. This framework can be iterated in various ways, for example, particular subsequences, rolling extrema, approximations of limits superior and inferior, the shape of the hills, the pre-image of $\|\color{black} G_N \color{black}\|_{\infty}$, and the relationship between support, variance, and kurtosis. 

We also presented computational evidence for a recurring pattern in the $\text{pdf}$: when $N=kq_n$ and $\frac{p_n}{q_n}$ is a good approximation for $\rho$, a predictable spiked-trapezoidal pattern emerges and eventually degrades. The precise onset of degradation and the threshold value of $\bigl|\rho - \tfrac{p_n}{q_n}\bigr|$ permitting persistence beyond $k=1$ \textcolor{black}{are left for future work}. Many other complex structures emerge as $N$ varies which can also be explored.

\subsection{Future work}
Numerically, the algorithm can be further refined. As noted in Section~\ref{sec:DTA}, the distances between consecutive \textcolor{black}{discontinuities} in $\color{black} G_N \color{black}(x,\rho)$ take on at most three \textcolor{black}{distinct} values - the third, if it exists, being the sum of the other two. This property could be exploited to improve both efficiency and accuracy.

The methodology developed in this paper can likely be extended to different types of sums,
\begin{align*}
   \color{black} S_N(x,\rho) &= \sum_{k=1}^N f\left( \{x+ k \rho \} \right) -\color{black} N\int_0^1 f(t) \,dt 
\end{align*}
\textcolor{black}{For irrational $\rho$, the sequence $\{x+ k \rho \}$ is uniformly distributed modulo one, so $\frac{S_N}{N} \rightarrow 0$ as $N \rightarrow \infty$ for every continuous function $f$~\cite{kuipers1974uniform}.
Of course, not all continuous functions are equally numerically tractable. A natural next step would be the quadratic polynomial $f(t)=\left(t-\frac{1}{2}\right)^2$, with $\int_0^1\left(t-\frac{1}{2}\right)^2 \,dt= \frac{1}{12} $. For $x_k = \{x+k\rho\}$, the resulting sum is}
\begin{align*}
   \color{black} S_N^{(2)}(x,\rho) &= \color{black} \sum_{k=1}^N \left(\{ x+k\rho\}-\frac{1}{2}\right)^2 - \frac{N}{12}\\
\end{align*}
\textcolor{black}{Mori et al. derived exact expressions and asymptotic estimates for $S_N^{(2)}(0,\rho)$~\cite{mori2019distribution} and studied geometric patterns over $N$. Using $f(t)=\left(t-\frac{1}{2}\right)^2$ rather than $f(t)=t^2$ eliminates jump discontinuities because $f(0)=f(1)=\frac{1}{4}$. Instead, $S_N^{(2)}(x,\rho)$ is continuous and piecewise quadratic. More generally, for polynomials with $f(0)=f(1)$, the resulting sum is continuous and piecewise polynomial, and a similar approach may be appropriate (tracking the piecewise polynomials rather than the jumps), although the $\text{pdf}$ and related statistical formulas would need to be re-derived. }

\section*{Acknowledgments}
\textcolor{black}{The author thanks Corey Reidy for careful proofreading, as well as assistance with grammar and language, and Christian Wei{\ss} for the introduction to much of the relevant mathematical context and literature, as well as for many helpful discussions. The author also thanks the anonymous referees for their careful reading and helpful suggestions, which improved the manuscript.}\\
\noindent
\textcolor{black}{Additionally, the author acknowledges the use of Google Gemini 3 Pro and OpenAI ChatGPT-5.6 for minor editing assistance and checking mathematical and bibliographic consistency during the preparation of this manuscript. The author reviewed all AI output and assumes full responsibility for the content of the publication.}

\section*{Code availability}
The MATLAB code used to generate the results in this paper is archived at Zenodo and is publicly available at DOI: \url{https://doi.org/10.5281/zenodo.17612959}. The repository contains all scripts and instructions necessary to reproduce the figures in the paper.

\section*{Disclosure statement}
The author reports there are no competing interests to declare.

\bibliographystyle{siamplain}
\bibliography{bibliography}

@article{birkhoff1931proof,
  title={Proof of the ergodic theorem},
  author={Birkhoff, George D},
  journal={Proceedings of the {N}ational {A}cademy of {S}ciences},
  volume={17},
  number={12},
  pages={656--660},
  year={1931},
  publisher={National {A}cademy of {S}ciences},
  doi={10.1073/pnas.17.2.656}
}

@article{knill2011self,
  title={Self-similarity and growth in {B}irkhoff sums for the golden rotation},
  author={Knill, Oliver and Tangerman, Folkert},
  journal={Nonlinearity},
  volume={24},
  number={11},
  year={2011},
  publisher={IOP Publishing},
  doi={10.1088/0951-7715/24/11/006}
}

@article{kochergin2023growth,
  title={On the Growth of {B}irkhoff Sums over a Rotation of the Circle},
  author={Kochergin, Andrey Vasilyevich},
  journal={Mathematical Notes},
  volume={113},
  number={5},
  pages={784--793},
  year={2023},
  publisher={Springer},
  doi={10.1134/S0001434623050206}
}

@article{cuyt2001remarkable,
  title={A remarkable example of catastrophic cancellation unraveled},
  author={Cuyt, Annie and Verdonk, Brigitte and Becuwe, Stefan and Kuterna, Peter},
  journal={Computing},
  volume={66},
  pages={309--320},
  year={2001},
  publisher={Springer},
  doi={10.1007/s006070170028}
}

@article{ramshaw1981discrepancy,
  title={On the discrepancy of the sequence formed by the multiples of an irrational number},
  author={Ramshaw, Lyle},
  journal={Journal of Number Theory},
  volume={13},
  number={2},
  pages={138--175},
  year={1981},
  publisher={Elsevier},
  doi={10.1016/0022-314X(81)90002-0}
}

@article{roccadas2011local,
  title={On the local discrepancy of (n$\alpha$)-sequences},
  author={Ro{\c{c}}adas, Lu{\'\i}s and Schoi{\ss}engeier, Johannes},
  journal={Journal of Number Theory},
  volume={131},
  number={8},
  pages={1492--1497},
  year={2011},
  publisher={Elsevier},
  doi={10.1016/j.jnt.2011.01.016}
}

@article{setokuchi2015discrepancy,
  title={On the discrepancy of irrational rotations with isolated large partial quotients: long term effects},
  author={Setokuchi, T},
  journal={Acta Mathematica {H}ungarica},
  volume={147},
  number={2},
  pages={368--385},
  year={2015},
  publisher={Springer},
  doi={10.1007/s10474-015-0563-0}
}

@article{bountis2020cauchy,
  title={Cauchy distributions for the integrable standard map},
  author={Bountis, Anastasios and Veerman, J J P and Vivaldi, Franco},
  journal={Physics Letters A},
  volume={384},
  number={26},
  pages={126659},
  year={2020},
  publisher={Elsevier},
  doi={10.1016/j.physleta.2020.126659}
}

@book{kuipers1974uniform,
  title={Uniform distribution of sequences},
  author={Kuipers, Lauwerens and Niederreiter, Harald},
  year={1974},
  publisher={John {W}iley \& {S}ons, {I}nc.},
  ISBN={0471510459}
}

@misc{veerman_BS,
	title = {Birkhoff Measures, {B}irkhoff Sums, and Discrepancies},
	author = {Ralston, D and Tangerman, F M and Veerman, J J P and Wu, H},
    journal={arXiv preprint arXiv:2511.22802},
    year={2025},
    url = {https://arxiv.org/pdf/2511.22802}
}

@book{pitman2012probability,
  title={Probability},
  author={Pitman, Jim},
  year={1993},
  publisher={New York: Springer-Verlag},
  ISBN={0387979743}
}

@article{mori2019distribution,
  title={On the distribution of partial sums of irrational rotations},
  author={Mori, Yoshiyuki and Shimaru, Naoto and Takashima, Keizo},
  journal={Periodica Mathematica Hungarica},
  volume={78},
  number={1},
  pages={88--97},
  year={2019},
  publisher={Springer},
  doi={10.1007/s10998-018-00273-y}
}

@article{shimaru2018discrepancies,
  title={On discrepancies of irrational rotations with several large partial quotients},
  author={Shimaru, N and Takashima, K},
  journal={Acta Mathematica Hungarica},
  volume={156},
  number={2},
  year={2018},
  pages = {449--458},
  doi={10.1007/s10474-018-0875-y}
}

@article{doi2017upper,
  title={An upper estimate for the discrepancy of irrational rotations},
  author={Doi, K and Shimaru, N and Takashima, K},
  journal={Acta Mathematica Hungarica},
  volume={152},
  number={1},
  pages={109--113},
  year={2017},
  publisher={Springer},
  doi={10.1007/s10474-017-0702-x}
}

@article{antonevich2022behaviour,
  title={Behaviour of {B}irkhoff sums generated by rotations of the circle},
  author={Antonevich, Anatolij Borisovich and Kochergin, Andrey Vasilyevich and Shukur, Ali Abdulhussein},
  journal={Matematicheskii Sbornik},
  volume={213},
  number={7},
  pages={3--38},
  year={2022},
  doi={10.4213/sm9356e}
}

@article{van1988three,
  title={The three gap theorem ({S}teinhaus conjecture)},
  author={Van Ravenstein, Tony},
  journal={Journal of the {A}ustralian Mathematical Society},
  volume={45},
  number={3},
  pages={360--370},
  year={1988},
  publisher={Cambridge University Press},
  doi={10.1017/S1446788700031062}
}

@article{sedgewick1978implementing,
  title={Implementing quicksort programs},
  author={Sedgewick, Robert},
  journal={Communications of the ACM},
  volume={21},
  number={10},
  pages={847--857},
  year={1978},
  publisher={ACM New York, NY, USA},
  doi={10.1145/359619.35963}
}

@book{heath2018scientific,
  title={Scientific Computing: An Introductory Survey, Revised Second Edition},
  author={Heath, Michael T},
  year={2018},
  publisher={SIAM},
  ISBN={9781611975574}
}

@inproceedings{shamos1976geometric,
  title={Geometric intersection problems},
  author={Shamos, Michael Ian and Hoey, Dan},
  booktitle={17th Annual Symposium on Foundations of Computer Science},
  pages={208--215},
  year={1976},
  organization={IEEE},
  doi={10.1109/SFCS.1976.16}
}

@article{shimaru2017discrepancies,
  title={On discrepancies of irrational rotations: an approach via rational rotation},
  author={Shimaru, Naoto and Takashima, Keizo},
  journal={Periodica Mathematica {H}ungarica},
  volume={75},
  number={1},
  pages={29--35},
  year={2017},
  publisher={Springer},
  doi={10.1007/s10998-016-0164-x}
}

@book{casella2024statistical,
  title={Statistical Inference},
  author={Casella, George and Berger, Roger},
  year={2024},
  edition={2},
  publisher={{C}hapman and {H}all/{CRC}},
  doi = {10.1201/9781003456285}
}

@article{setokuchi2014discrepancies,
  title={Discrepancies of irrational rotations with isolated large partial quotients},
  author={Setokuchi, T and Takashima, K},
  journal={Uniform Distribution Theory},
  volume={9},
  number={2},
  pages={31--57},
  year={2014},
  note        = {URL: \url{https://pcwww.liv.ac.uk/~karpenk/JournalUDT/vol09/no2/03SetokuchiTakashima.pdf}}

}

@techreport{McKeeman2004Sort,
  author      = {McKeeman, Bill and Shure, Loren},
  title       = {An Adventure of Sorts--Behind the Scenes of a {MATLAB} Upgrade},
  institution = {{MathWorks}},
  year        = {2004},
  month       = {Oct},
}

@article{schoissengeier1986discrepancy,
  title={On the discrepancy of (n$\alpha$), {II}},
  author={Schoissengeier, Johannes},
  journal={Journal of Number Theory},
  volume={24},
  number={1},
  pages={54--64},
  year={1986},
  publisher={Elsevier},
  doi={10.1016/0022-314X(86)90057-0}
}

@software{MATLAB:R2025b,
  author       = {{The MathWorks, Inc.}},
  title        = {{MATLAB} version {R}2025b},
  year         = {2025},                         
  publisher    = {{The MathWorks, Inc.}},
  address      = {Natick, Massachusetts, United States},
  url          = {https://www.mathworks.com}
}

@article{blessing2024weighted,
  title={Weighted {B}irkhoff averages and the parameterization method},
  author={Blessing, David and Mireles James, JD},
  journal={{SIAM} Journal on Applied Dynamical Systems},
  volume={23},
  number={3},
  pages={1766--1804},
  year={2024},
  publisher={{SIAM}},
  doi={10.1137/23M1579546}
}

@article{gonzalez2022efficient,
  title={Efficient and reliable algorithms for the computation of non-twist invariant circles},
  author={Gonz{\'a}lez, Alejandra and Haro, {\`A}lex and de la Llave, Rafael},
  journal={Foundations of Computational Mathematics},
  volume={22},
  pages={791--847},
  year={2022},
  publisher={Springer},
  doi={10.1007/s10208-021-09517-9}
}

@incollection{das2016quasiperiodicity,
  title={Quasiperiodicity: rotation numbers},
  author={Das, Suddhasattwa and Saiki, Yoshitaka and Sander, Evelyn and Yorke, James A},
  booktitle={The Foundations of Chaos Revisited: From {P}oincar{\'e} to Recent Advancements},
  pages={103--118},
  year={2016},
  publisher={Springer},
  doi={10.1007/978-3-319-29701-9_7}
}

@article{fukuyama2016metric,
  title={Metric discrepancy results for geometric progressions with large ratios},
  author={Fukuyama, Katusi and Yamashita, Mai},
  journal={Monatshefte f{\"u}r Mathematik},
  volume={180},
  number={4},
  pages={731--742},
  year={2016},
  publisher={Springer},
  doi={10.1007/s00605-015-0791-y}
}

@article{fukuyama2023metric,
  title={Metric discrepancy results for geometric progressions with large ratios {II}},
  author={Fukuyama, Katusi},
  journal={Monatshefte f{\"u}r Mathematik},
  volume={202},
  number={2},
  pages={281--315},
  year={2023},
  publisher={Springer},
  doi={10.1007/s00605-023-01823-4}
}

@book{Khinchin_CF,
  title={Continued Fractions},
  author={Khinchin, A. Ya.},
  edition={3},
  year={1964},
  publisher={University of Chicago Press},
  address={Chicago},
  isbn={0226447499}
}

@article{vinson2001partial,
  title={Partial Sums of $\zeta$ ($1/2$) Modulo 1},
  author={Vinson, Jade},
  journal={Experimental Mathematics},
  volume={10},
  number={3},
  pages={337--344},
  year={2001},
  publisher={Taylor \& Francis},
  doi = {10.1080/10586458.2001.10504454}
}

@article{tong2024exponential,
  title={Exponential convergence of the weighted {B}irkhoff average},
  author={Tong, Zhicheng and Li, Yong},
  journal={Journal de Math{\'e}matiques Pures et Appliqu{\'e}es},
  volume={188},
  pages={470--492},
  year={2024},
  publisher={Elsevier},
  doi={10.1016/j.matpur.2024.06.003}
}

@article{dolgopyat2017temporal,
  title={Temporal distributional limit theorems for dynamical systems},
  author={Dolgopyat, Dmitry and Sarig, Omri},
  journal={Journal of Statistical Physics},
  volume={166},
  number={3},
  pages={680--713},
  year={2017},
  publisher={Springer},
  doi={10.1007/s10955-016-1689-3}
}

@article{dolgopyat2020quenched,
  title={Quenched and annealed temporal limit theorems for circle rotations},
  author={Dolgopyat, Dmitry and Sarig, Omri},
  journal={Ast\'erisque},
  volume={415},
  number={415},
  pages={59--85},
  year={2020},
  publisher={Societe Mathematique de France},
  doi={10.24033/ast.1100}
}

@article{bromberg2018temporal,
  title={A temporal central limit theorem for real-valued cocycles over rotations},
  author={Bromberg, Michael and Ulcigrai, Corinna},
  journal={Annales de l'Institut Henri Poincar{\'e} (B) Probabilit{\'e}s et Statistiques},
  volume={54},
  number={4},
  pages={2304--2334},
  year={2018},
  doi={10.1214/17-AIHP872}
}

@article{fruhwirth2024birkhoff,
  title={On {B}irkhoff sums that satisfy no temporal distributional limit theorem for almost every irrational},
  author={Fr{\"u}hwirth, Lorenz and Hauke, Manuel},
  journal={Annales Henri Lebesgue},
  volume={7},
  pages={251--265},
  year={2024},
  doi={10.5802/ahl.200}
}

@article{dolgopyat2018no,
  title={No temporal distributional limit theorem for a.e. irrational translation},
  author={Dolgopyat, Dmitry and Sarig, Omri},
  journal={Annales {H}enri {L}ebesgue},
  volume={1},
  pages={127--148},
  year={2018},
  doi={10.5802/ahl.4}
}

@article{das2018super,
  title={Super convergence of ergodic averages for quasiperiodic orbits},
  author={Das, Suddhasattwa and Yorke, James A},
  journal={Nonlinearity},
  volume={31},
  number={2},
  pages={491--501},
  year={2018},
  publisher={IOP Publishing},
  doi={10.1088/1361-6544/aa99a0}
}

@article{tong2026weighted,
  title={Weighted {B}irkhoff averages: Deterministic and probabilistic perspectives},
  author={Tong, Zhicheng and Li, Yong},
  journal={Mathematische Annalen},
  volume={394},
  number={1},
  pages={18},
  year={2026},
  publisher={Springer},
  doi={10.1007/s00208-026-03311-0}
}

@article{pinner1997sums,
  title={On sums of fractional parts $\{$n$\alpha$+ $\gamma$$\}$},
  author={Pinner, Christopher},
  journal={Journal of Number Theory},
  volume={65},
  number={1},
  pages={48--73},
  year={1997},
  publisher={Elsevier},
  doi={10.1006/jnth.1997.2080}
}

@article{krasil2011description,
  title={Description and exact maximum and minimum values of the remainder in the problem of the distribution of fractional parts},
  author={Krasil’shchikov, Vasilij Vyacheslavovich and Shutov, Anton Vladimirovich},
  journal={Mathematical Notes},
  volume={89},
  number={1},
  pages={59--67},
  year={2011},
  publisher={Springer},
  doi={10.1134/S0001434611010068}
}

@article{baxa1994minimum,
  title={Minimum and maximum order of magnitude of the discrepancy of (n$\alpha$)},
  author={Baxa, C and Schoissengeier, J},
  journal={Acta Arithmetica},
  volume={68},
  number={3},
  pages={281--290},
  year={1994},
  doi={10.4064/aa-68-3-281-290}
}

@article{borda2023limit,
  title={Limit laws of maximal {B}irkhoff sums for circle rotations via quantum modular forms},
  author={Borda, Bence},
  journal={International Mathematics Research Notices},
  volume={2023},
  number={22},
  pages={19340--19389},
  year={2023},
  publisher={Oxford University Press},
  doi={10.1093/imrn/rnad107}
}

@article{niederreiter1978quasi,
  title={Quasi-{M}onte {C}arlo methods and pseudo-random numbers},
  author={Niederreiter, Harald},
  journal ={Bulletin of the American Mathematical Society},
  volume={84},
  number={6},
  year={1978},
  pages = {957--1041},
  doi={10.1090/S0002-9904-1978-14532-7}
}

@article{polanco2016continuous,
  title={Continuous distributions arising from the three gap theorem},
  author={Polanco, Gerem{\'\i}as and Schultz, Daniel and Zaharescu, Alexandru},
  journal={International Journal of Number Theory},
  volume={12},
  number={07},
  pages={1743--1764},
  year={2016},
  publisher={World Scientific},
  doi={10.1142/S1793042116501074}
}

@article{ravenstein1990characteristics,
  title={Characteristics and the three gap theorem},
  author={Ravenstein, Tony van and Winley, Graham and Tognetti, Keith},
  journal={The Fibonacci Quarterly},
  volume={28},
  number={3},
  pages={204--214},
  year={1990},
  publisher={Taylor \& Francis},
  doi={10.1080/00150517.1990.12429482}
}

\newpage 
\appendix

\section{Selected statistical plots} \label{apx:stats_plots}
Statistical plots are shown for $\rho$ of various \textcolor{black}{arithmetic} types for $N=2, 3, 4 \dots$.

\begin{figure}[H]
    \centering
    \includegraphics[angle=-90, width=0.8\linewidth]{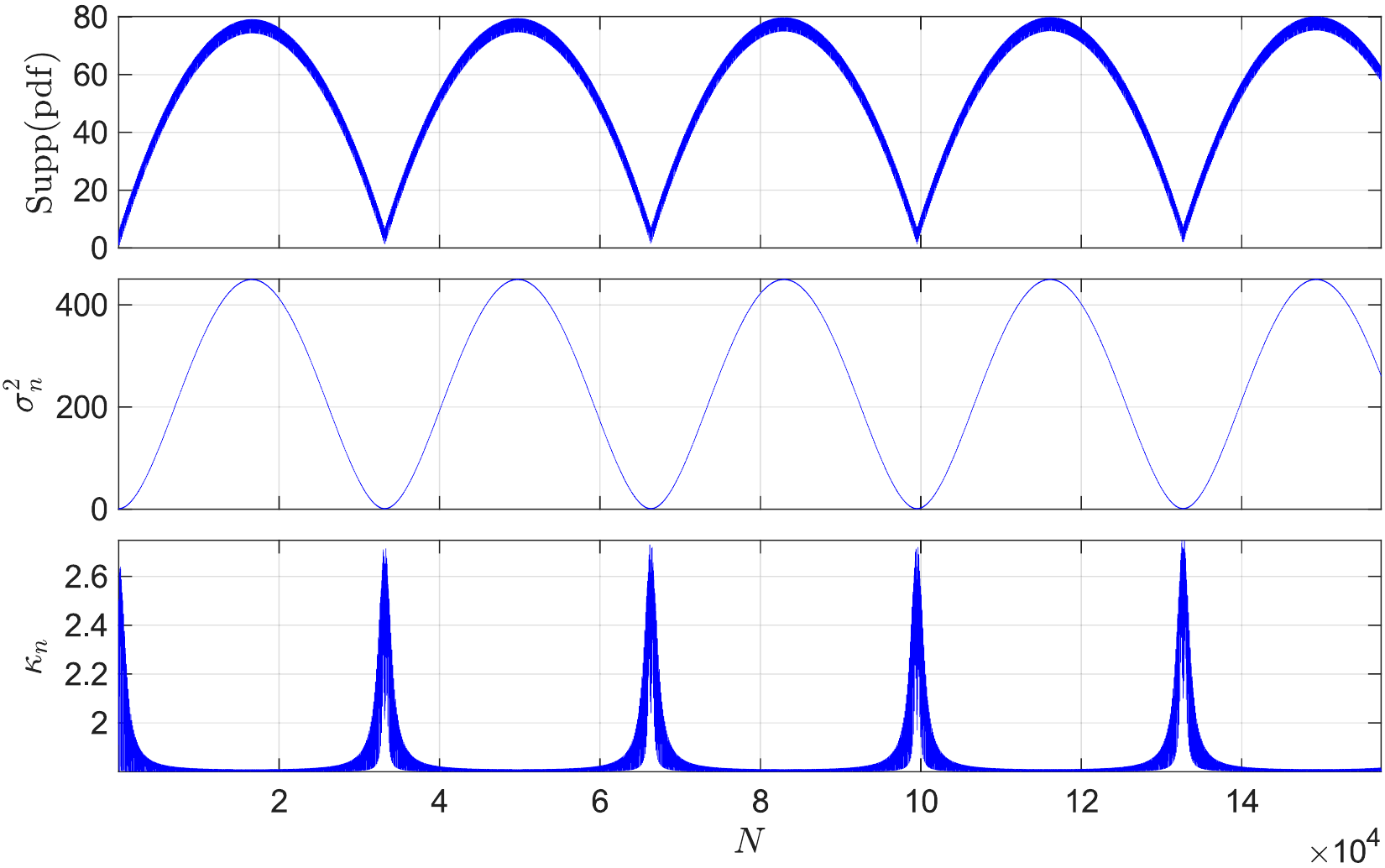}
    \caption{Statistical plots for $\rho=\pi-3$} \label{fig:pi_stats}
\end{figure}

\begin{figure}[H]
    \centering
    \includegraphics[angle=-90,width=0.8\linewidth]{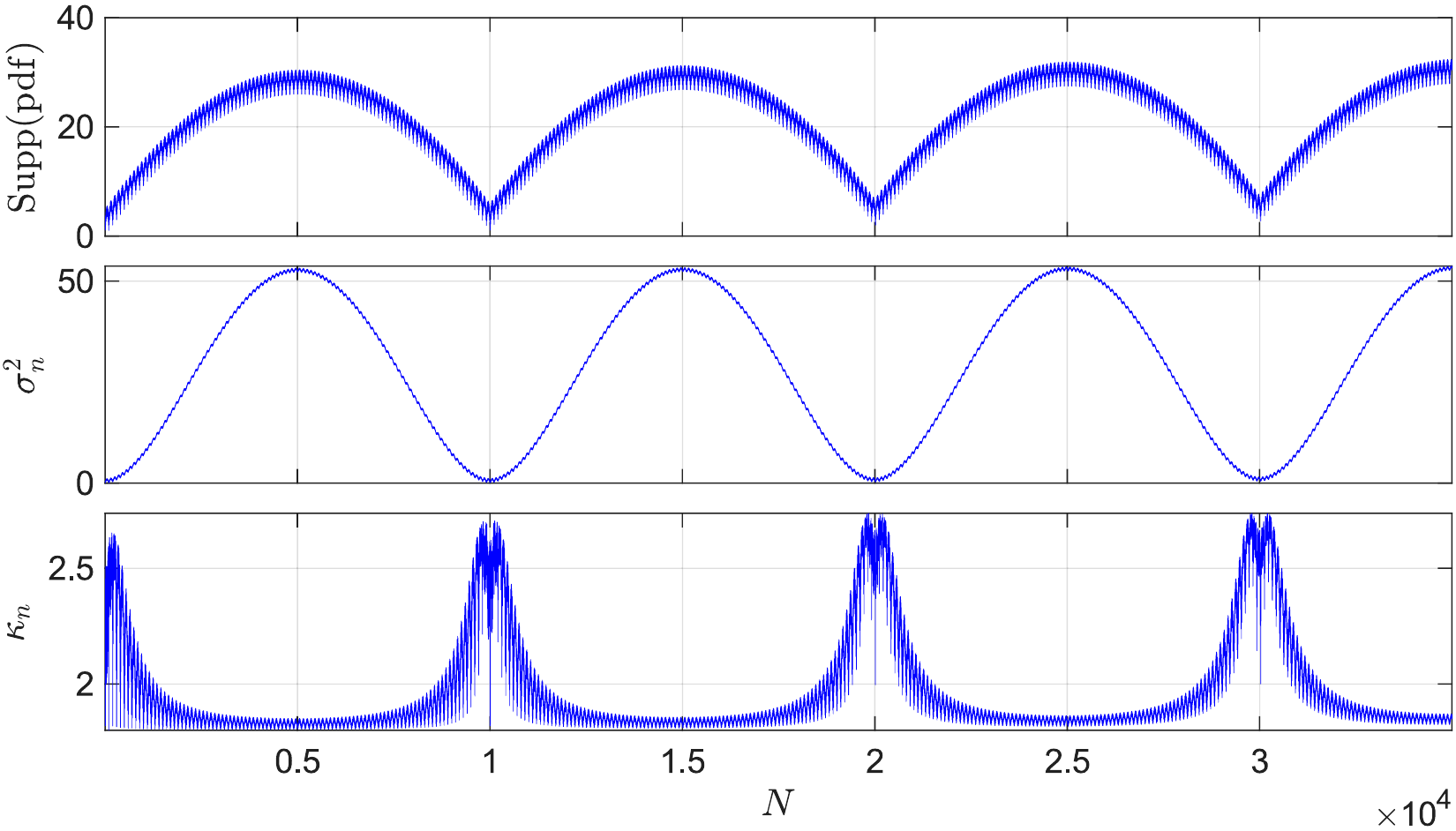}
    \caption{Statistical plots for $\rho=L$} \label{fig:lio_stats}
\end{figure}


\begin{figure}[H]
    \centering
    \includegraphics[angle=-90, width=0.8\linewidth]{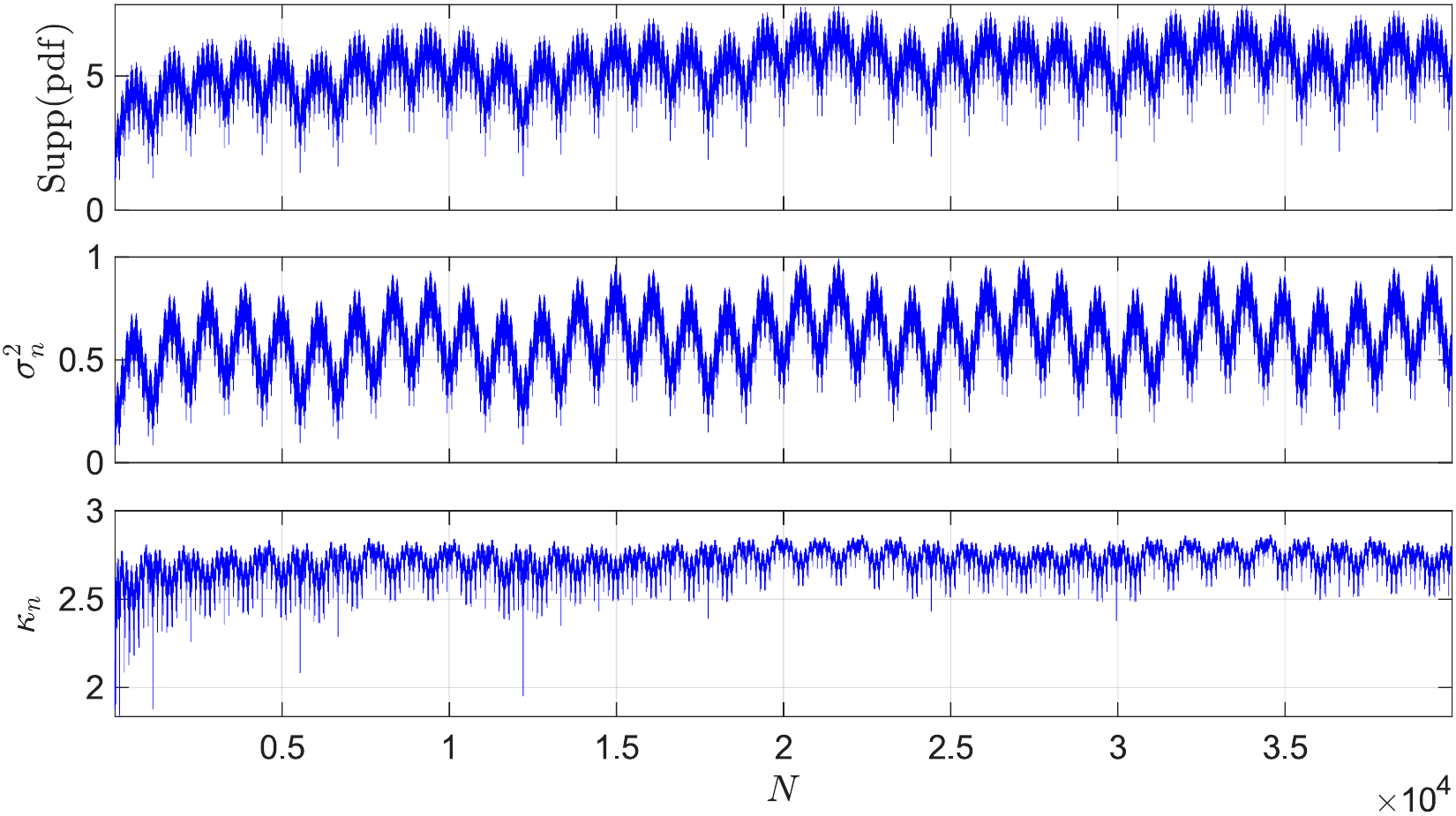}
    \caption{Statistical plots for $\rho=\frac{\pi^2}{6}-1$} \label{fig:zeta_stats}
\end{figure}

\begin{figure}[H]
    \centering
    \includegraphics[angle=-90, width=0.8\linewidth]{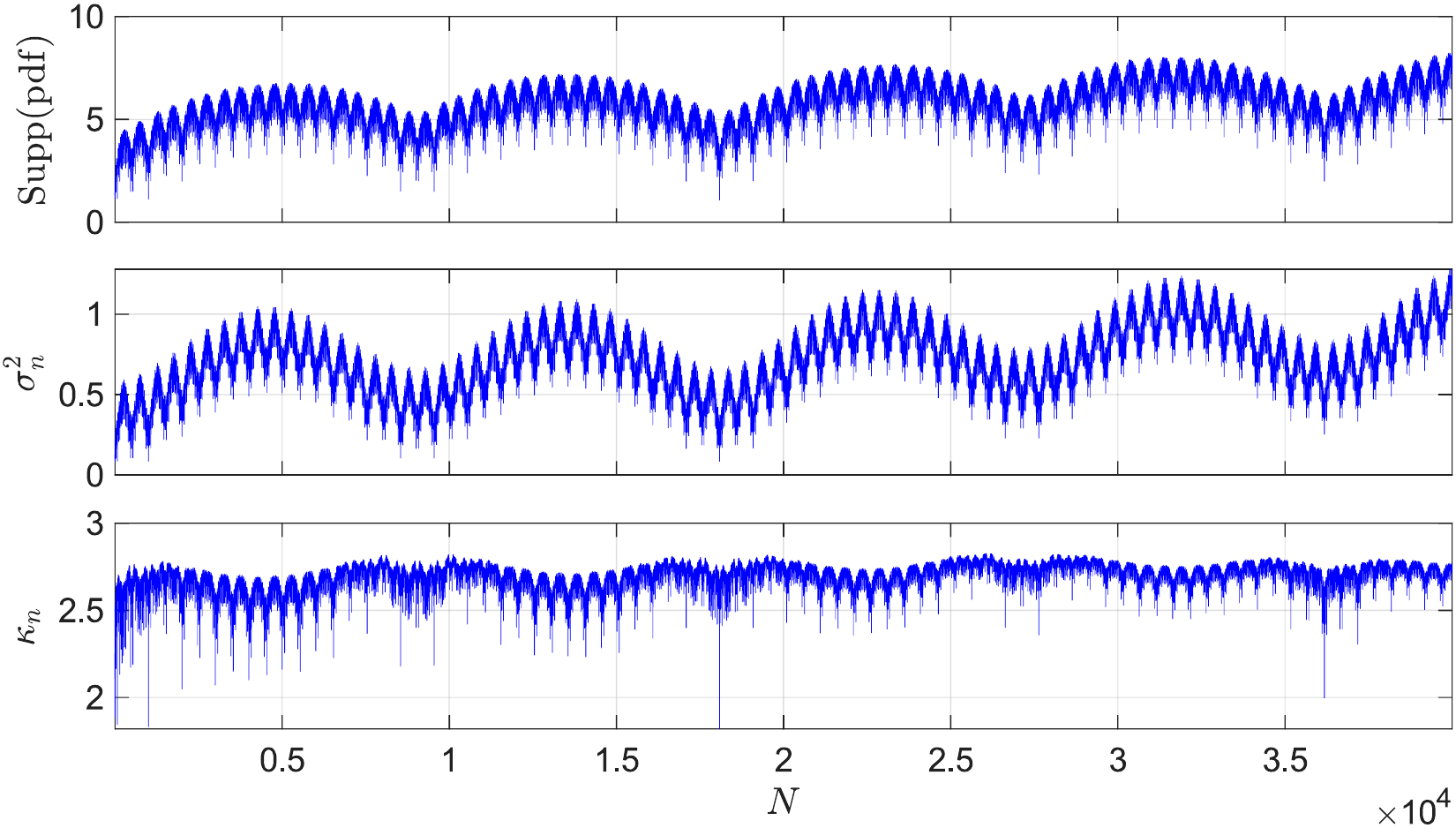}
    \caption{Statistical plots for $\rho=e-2$} \label{fig:exp_stats}
\end{figure}


\begin{figure}[H]
    \centering
    \includegraphics[angle=-90, width=0.8\linewidth]{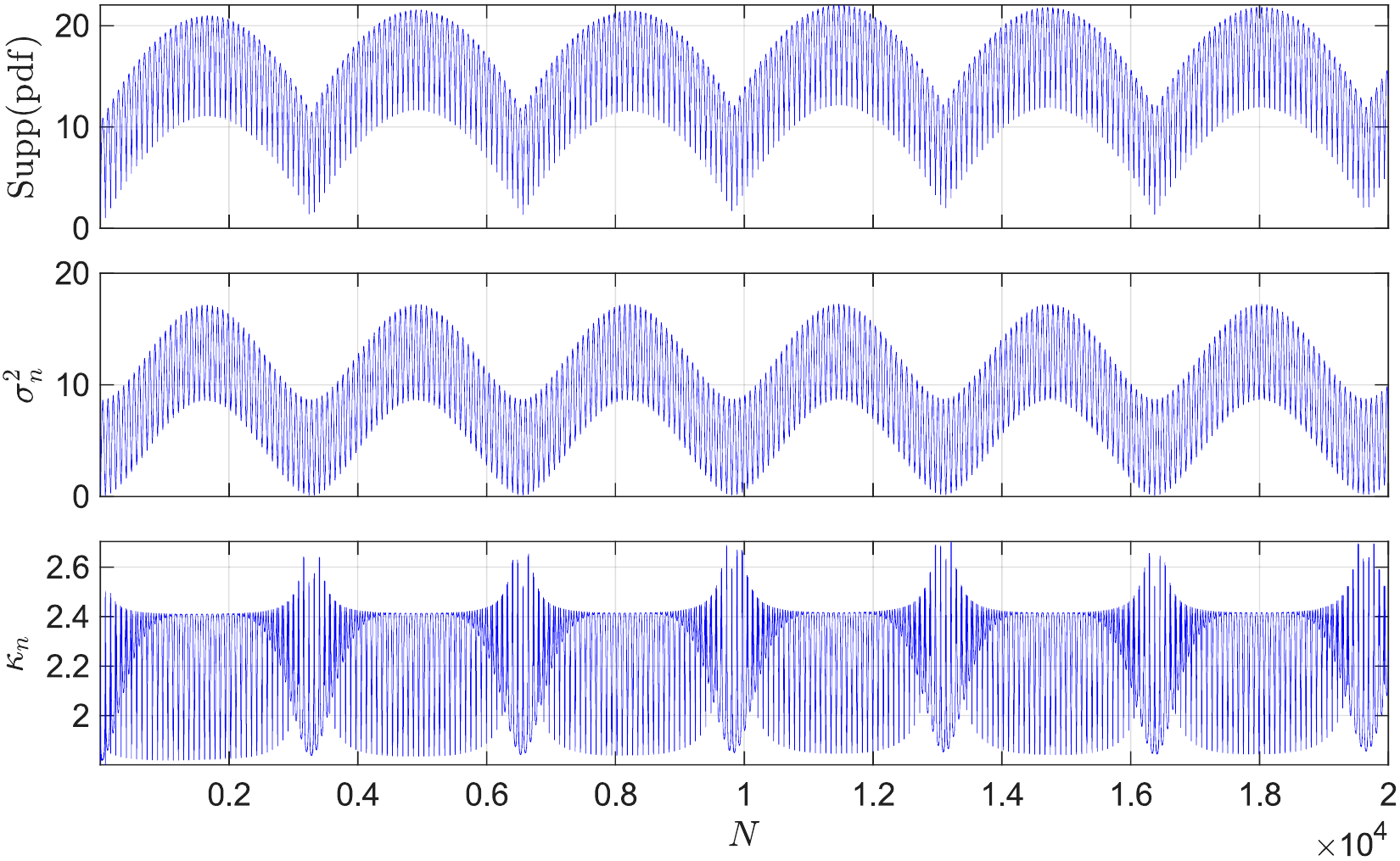}
    \caption{Statistical plots for $\rho=\lbrack \color{black} 0; \color{black}2, 40, 40, \color{black} \overline{2} \color{black} \rbrack$} \label{fig:rho_2_40402_stats}
\end{figure}

\begin{figure}[H]
    \centering
    \includegraphics[angle=-90, width=0.8\linewidth]{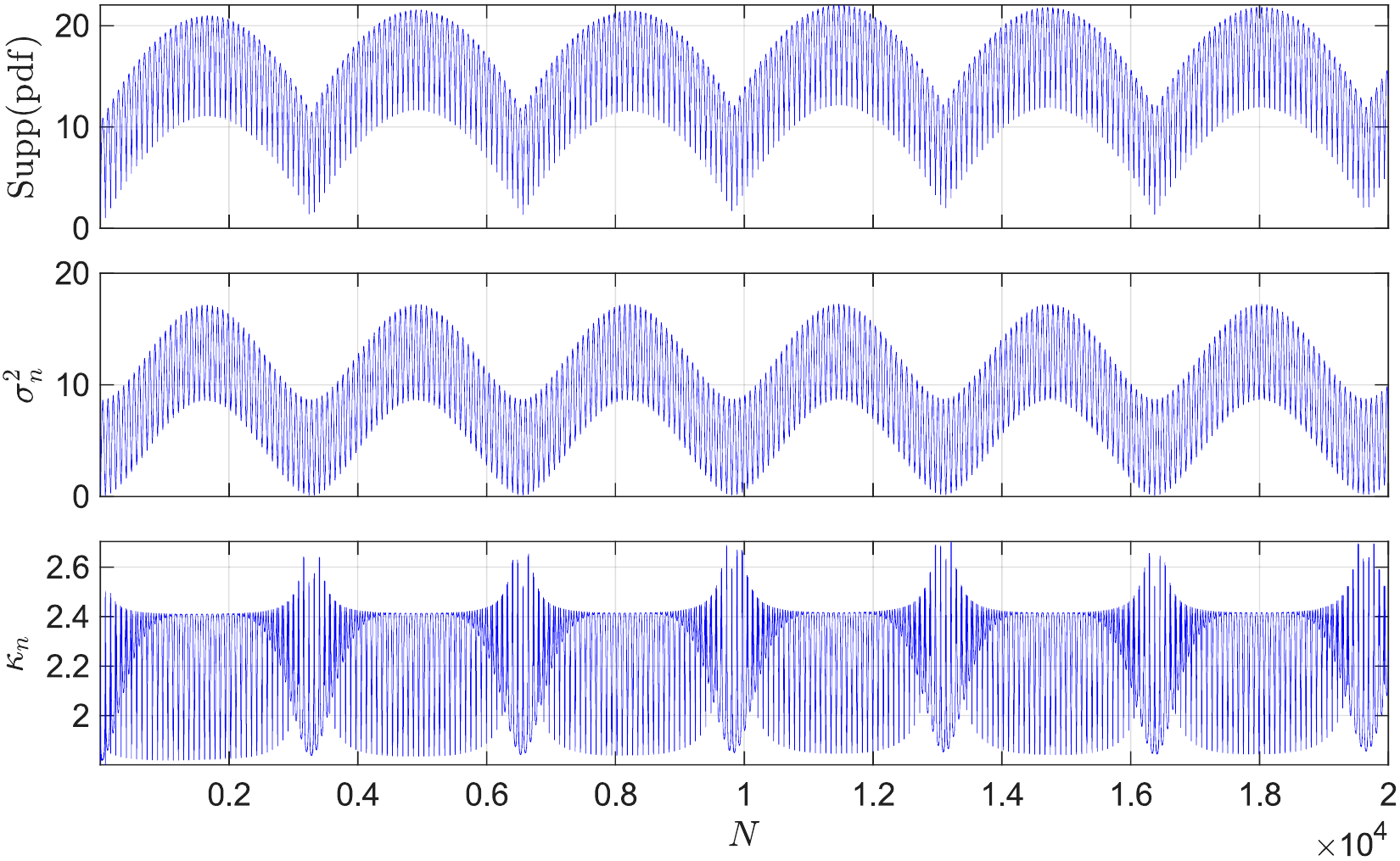}
    \caption{Statistical plots for $\rho=\lbrack\color{black} 0; \color{black} 2, 2, 2, 2, 10, \color{black} \overline{2} \color{black} \rbrack$} \label{fig:rho_2_2210_stats}
\end{figure}

\begin{figure}[H]
    \centering
    \includegraphics[angle=-90, width=0.8\linewidth]{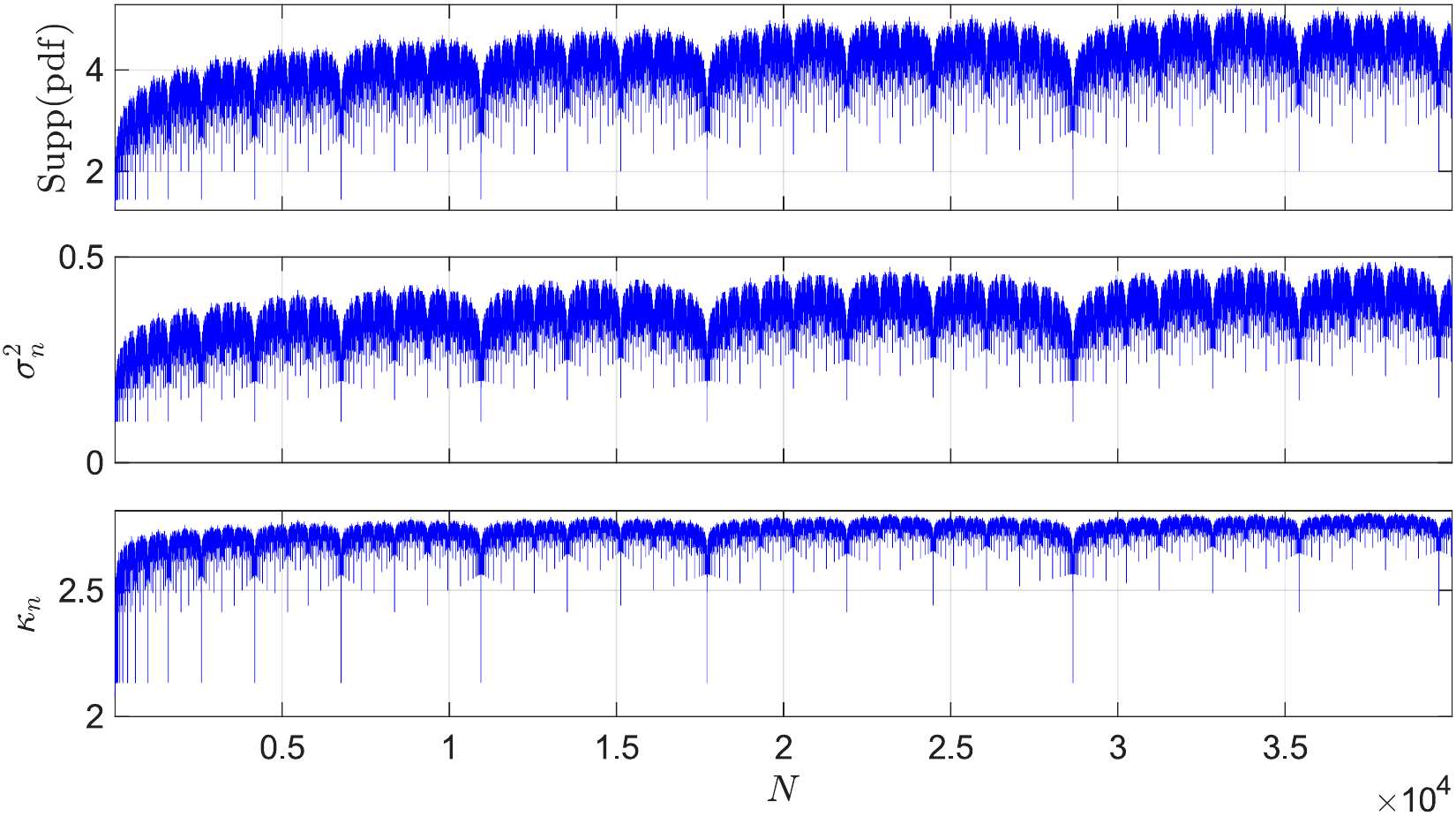}
    \caption{Statistical \textcolor{black}{plots} for $\rho=\frac{\sqrt{5}-1}{2}$} \label{fig:golden_stats}
\end{figure}

\begin{figure}[H]
    \centering
    \includegraphics[angle=-90, width=0.8\linewidth]{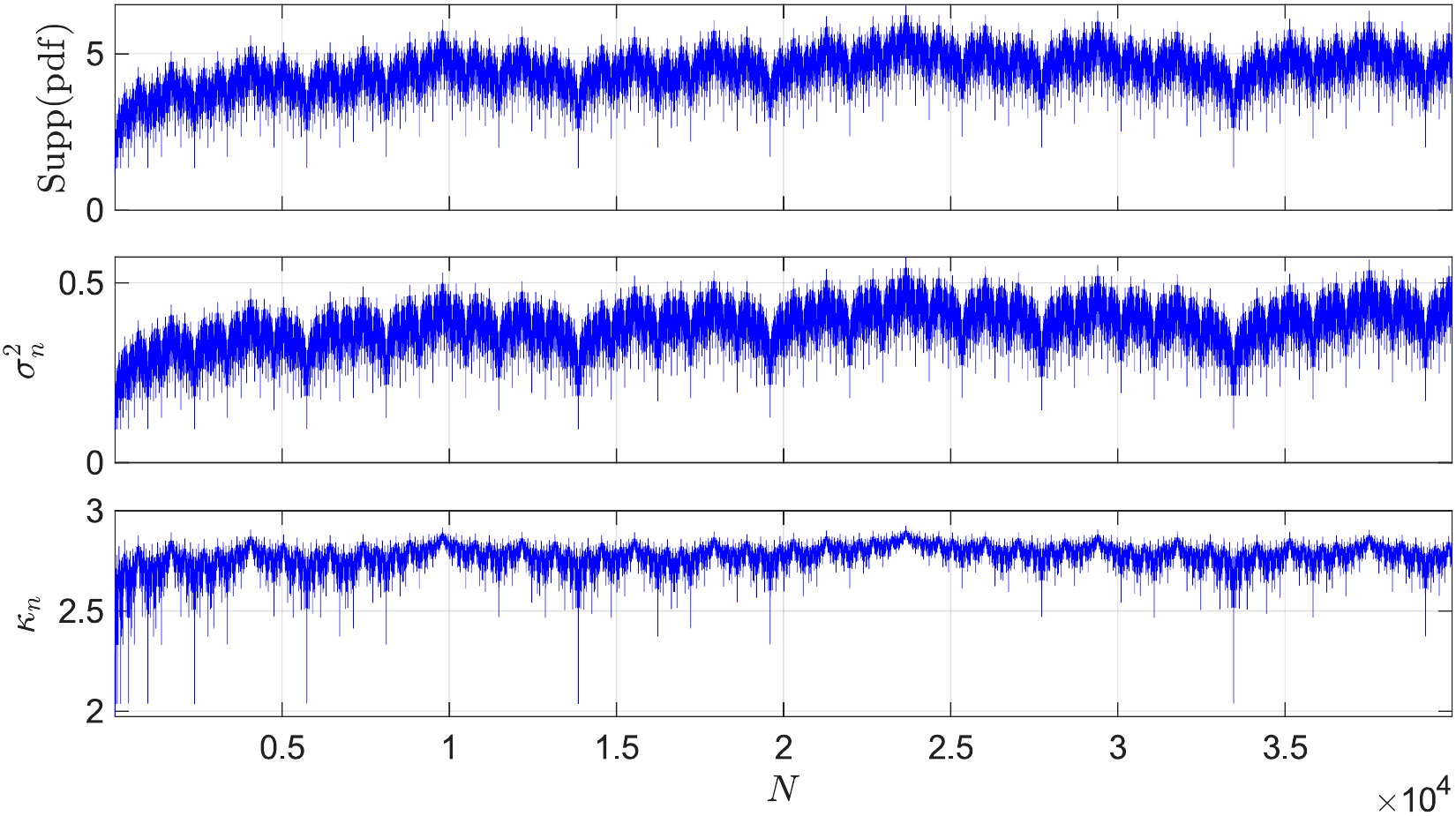}
    \caption{Statistical plots for $\rho=\sqrt{2}-1$} \label{fig:silver_stats}
\end{figure}

\end{document}